\documentclass{jnmp}

\usepackage{amsmath}
\usepackage{lscape}

\JNMPnumberwithin{equation}{section}

\DeclareMathOperator{\Hom}{Hom}
\DeclareMathOperator{\Lie}{Lie}
\DeclareMathOperator{\id}{id}
\DeclareMathOperator{\Vect}{Vect}
\DeclareMathOperator{\diag}{diag}
\DeclareMathOperator{\Span}{Span}
\DeclareMathOperator{\ad}{ad}
\DeclareMathOperator{\End}{End}
\DeclareMathOperator{\Par}{Par}
\DeclareMathOperator{\tr}{tr}
\DeclareMathOperator{\str}{str}

\DeclareMathOperator{\ber}{ber}
\DeclareMathOperator{\Vol}{Vol}

\DeclareMathOperator{\Bil}{Bil}
\DeclareMathOperator{\antidiag}{antidiag}
\DeclareMathOperator{\Spec}{Spec}

\newtheorem{theorem}{Theorem}
\newtheorem{proposition}{Proposition}
\newtheorem{lemma}{Lemma}
\newtheorem{corollary}{Corollary}

\theoremstyle{definition}

\newtheorem{remark}{Remark}
\newtheorem*{example}{Example}

\begin{document}

\renewcommand{\evenhead}{D~Leites, E~Poletaeva and V~Serganova}
\renewcommand{\oddhead}{On Einstein Equations on Manifolds and Supermanifolds}

\thispagestyle{empty}

\FirstPageHead{9}{4}{2002}{\pageref{leites-firstpage}--\pageref{leites-lastpage}}{Article}

\copyrightnote{2002}{D~Leites, E~Poletaeva and V~Serganova}

\Name{On Einstein Equations on Manifolds\\
and Supermanifolds}

\label{leites-firstpage}

\Author{D~LEITES~${}^{\dag^1}$, E~POLETAEVA~${}^{\dag^2}$ and V~SERGANOVA~${}^{\dag^3}$}

\Address{${}^{\dag^1}$~Department of Mathematics, University
of Stockholm, Roslagsv.  101, \\
$\phantom{{}^{\dag^1}}$~Kr\"aftriket hus 6, SE-106 91, Stockholm, Sweden\\
$\phantom{{}^{\dag^1}}$~E-mail: mleites@matematik.su.se\\[10pt]
${}^{\dag^2}$~Department of Mathematics, Lund University, Sweden\\
 $\phantom{{}^{\dag^2}}$~E-mail: elena@maths.lth.se\\[10pt]
${}^{\dag^3}$~Department of Mathematics, University of California, Berkeley, USA\\
$\phantom{{}^{\dag^3}}$~E-mail: serganov@math.berkeley.edu}

\Date{Received February 27, 2001; Revised August 29, 2001;
Accepted April 11, 2002}

\begin{abstract}
\noindent The Einstein equations (EE) are certain conditions on the
Riemann tensor on the real Minkowski space $M$.  In the twistor
picture, after complexification and compactifi\-cation~$M$ becomes the
Grassmannian $Gr_{2}^{4}$ of $2$-dimensional subspaces in the
$4$-di\-men\-sional complex one.  Here we answer for which of the
classical domains considered as manifolds with $G$-structure it is
possible to impose conditions similar in some sense to~EE. The above
investigation has its counterpart on superdomains: an analog of the
Riemann tensor is defined for any {\it super}manifold with
$G$-structure with any Lie supergroup~$G$.  We also derive similar
analogues of EE on supermanifolds.  Our analogs of EE are not what
physicists consider as SUGRA (supergravity), for SUGRA see~\cite{GL4,LP2}.
\end{abstract}

\section{Introduction}

This is an expanded version of a part of Leites' lectures at ICTP,
Trieste, in March 1991 on our results.  The description of ``the left
hand side of $N$-extended SUGRA equations'', though computed several
years later, appeared earlier~\cite{GL1} and refers to some results
from this paper and~\cite{LP2}.

Roughly speaking, in this paper, as well as in \cite{GL1,P1,P2,P3,P4},
for a ${\mathbb Z}$-graded Lie superalgebra ${\mathfrak
g}_{*}=\mathop{\oplus}\limits_{i \geq -d}{\mathfrak g}_i$ and its
subalgebra ${\mathfrak g}_{-}=\mathop{\oplus}\limits_{i<
0}{\mathfrak g}_i$ we calculate $H^k({\mathfrak g}_{-}; {\mathfrak
g}_{*})$ for $k\leq 2$. In addition to a new result (analogs of EE)
this paper contains a summary of \cite{P1,P2,P3}.
The Nijenhuis tensor deserves a separate publication \cite{P4}.

For ${\mathfrak g}_{*}$ simple, $k=2$ and $d=1$, this cohomology can
be interpreted as analogs of the conformal part of the Riemann tensor,
called the {\it Weyl tensor}, (more exactly, values thereof at
a~point).  This cohomology coincides with the Weyl tensor on the
$n$-dimensional manifold when ${\mathfrak g}_{*}={\mathfrak o}(n+2)$.
For other Lie algebras ${\mathfrak g}_{*}$, not only simple ones, we
obtain the so-called {\it structure functions} (obstructions to
non-flatness in some sense) of the manifold with the $G$-structure,
where ${\mathfrak g}_{0}=\Lie (G)$, the Lie algebra of $G$.  For
${\mathfrak g}_{*}$ simple, the cohomology $H^k({\mathfrak g}_{-1};
{\mathfrak g}_{-1}\oplus \hat {\mathfrak g}_{0})$, where $\hat
{\mathfrak g}_{0}$ is the semisimple part of ${\mathfrak g}_{0}$,
corresponds to (the analogs of) the {\it Riemann tensor}; they consist
of $H^k({\mathfrak g}_{-}; {\mathfrak g}_{*})$ --- the ``conformal
part'' --- plus something else, and it is this ``extra'' part that
plays the main role in the left hand side of the Einstein equations.

For $d>1$ one obtains new invariants which we interpret as
obstructions to ``non-flatness'' of a manifold (or supermanifolds)
with a nonholonomic structure, see \cite{L7,LP2,GL1,GL2}.  These
invariants eluded researchers for almost a century, see Vershik's
review \cite{V}, where doubted if they existed. Similar structures
appear in Manin's book~\cite{M}, and our approach shows a method to
describe their ``non-flatness''. We have only started to study such structures; the detailed exposition is in preparation.

In this paper we only consider $d=1$ and mostly finite dimensional
cases.  Goncharov considered Lie algebras and cases when~${\mathfrak
g}_{*}$ is simple; we consider also superalgebras for ${\mathfrak
g}_{*}$ simple or close to simple and also consider $H^k({\mathfrak
g}_{-1}; {\mathfrak g}_{-1}\oplus \hat {\mathfrak g}_{0})$.  Other
cases are either open problems or will be considered elsewhere.

We are thankful to Grozman who verified our calculations of structure
functions for the exceptional superdomains and in several other cases
by means of his SuperLie package: even these finite dimensional
calculations are almost impossible to perform without computer whereas
to Grozman's package this is a matter of minutes in components; to
glue the components into a module takes several hours in each case.

The main object in the study of Riemannian geometry is the {\it
Riemann tensor}.  Under the action of $O(n)$ the space of values of
the Riemann tensor at the point splits into irreducible components
called the {\it Weyl tensor},  the {\it traceless Ricci tensor} and the
{\it scalar curvature}.  (On 4-dimensional manifolds the Weyl tensor
additionally splits into 2 subcomponents.)  All these tensors are
obstructions to the possibility of ``flattening" the canonical
(Levi--Civita) connection on the manifold they are considered.

More generally, let $G$ be any Lie group, not necessarily $O(n)$.  In
what follows we will recall the definition of a $G$-structure on a
manifold $M$ and structure functions of this $G$-structure.  Structure
functions are obstructions to integrability or, in other words, to the
possibility of ``flattening'' the $G$-structure or a connection
associated with it, sometimes, canonically, see~\cite{G}.  The Riemann
tensor is the only nontrivial
structure function for $G=O(n)$.  Several most known (or
popular recently) examples of $G$-structures and respective tensors
are:
\begin{center}\small
\begin{tabular}{|c|c|c|}
\hline
Name of the structure&$G$&Name of the tensor\\
\hline\hline
almost conformal&$G = CO(n)= O(n)\times {\mathbb R} ^+$&Weyl  tensor\\
\hline
Riemannian structure&$G = O(n)$&Riemann tensor\\
\hline
Penrose' twistors&$G = S(GL(2; \mathbb C)\times GL(2; \mathbb C))$&the
``$\alpha$-forms" and ``$\beta$-forms"\\
\hline
almost complex structure&$G = GL(n;{\mathbb C})\subset GL(2n;{\mathbb R} )$&Nijenhuis  tensor\\
\hline
almost symplectic structure&$G = Sp(2n)$&no accepted name\\
\hline
\end{tabular}
\end{center}

\begin{remark} The adverb ``almost" should be always added
until the $G$-structure under study is proved flat,  i.e.,  integrable;
by abuse of language people often omit it.
\end{remark}

Infinitesimal automorphisms (with polynomial coefficients)
of the flat $G$-structure on ${\mathbb R}^n$ ($n = \dim M$)
constitute the Cartan prolong
(see Section~2.2) --- the Lie algebra
$({\mathfrak g}_{-1}, {\mathfrak g}_{0})_{*}$, where ${\mathfrak
g}_{-1}$ can be identified with the tangent space $T_{m}M$ at a point
and ${\mathfrak g}_{0}=\Lie (G)$.  We interpret structure functions as
certain Lie algebra cohomology associated with $({\mathfrak g}_{-1},
{\mathfrak g}_{0})_{*}$.

The Riemannian case is the {\it reduction of the structure group} of
the conformal case.  More generally, if ${\mathfrak g}_{0}$ is a
central extension of the (semi-)simple Lie algebra $\hat{\mathfrak
g}_{0}$, the corresponding structure functions will be called, after
Goncharov, generalized conformal ones, whereas the structure functions
for $\hat{\mathfrak g}_{0}$ will correspond to a generalized
Riemannian case~--- a~possible candidate in search for analogs of
Einstein equations.

In \cite{G} Goncharov calculated all structure functions for the
analogues of the almost conformal structure corresponding to
irreducible compact Hermitian symmetric spaces (shortly CHSS in what
follows); in Goncharov's examples $G$ is the reductive part of the
stabilizer of any point of the compact Hermitian symmetric space.

Here we will consider the reductions of the cases considered by
Goncharov~--- analogues of Riemannian structures and their various
generalizations to ma\-ni\-folds and super\-ma\-ni\-folds, in
particular, infinite dimensional, associated with Kac--Moody and
stringy (super)algebras.  We also review some cases considered in
\cite{P1} and give an overview of \cite{P2} and \cite{LP1,LP2}.

For prerequisites on symmetric spaces see~\cite{H}.  Appendix contains
preliminaries on Lie superalgebras and supermanifolds; the
super analogs of classical symmetric spaces listed in \cite{S,LSV1}
are recalled in Tables.  Observe several interesting points.

(1) Some of the spaces and superspaces we distinguish are infinite
dimensional.  Some of these infinite dimensional analogues of EE can
only be realized on the total spaces of Fock bundles over
supermanifolds with at least~3 odd coordinates; the invariance group
of such an EE contains a contact Lie superalgebra.  Other infinite
dimensional examples are associated with Kac--Moody or loop algebras
of which the examples associated with twisted versions are most
intriguing.

(2) On supermanifolds, our analogues of EE are {\it not} what
physicists consider as {\it supergravity equations} (SUGRA); each
$N$-extended SUGRA requires a nonholonomic distribution and they are
considered in \cite{LP2,GL1}. Recall that having struggled for a
decade with a~conventional model of Minkowski superspace for deriving
$N=2$ SUGRA the Ogievetsky's group GIKOS had found a
solution~\cite{GI}: one has to enlarge the Minkowski space underlying
Minkowski superspace for $N=2$ with an additional ``harmonic'' space
$P^1$.  How to advance as $N$ grows was unclear, cf.\ pessimistic
remarks in \cite{GI} and \cite{WZ}.

What was the problem?!  Take the usual recipe for calculation of the
Riemann tensor or even structure function of any $G$-structure,
insert some signs to account for super flavour and that will be it!
This is more or less what is suggested in \cite{CDF} and \cite{RS}.
The snag is that in doing so we tacitly assume that ${\mathfrak
g}_{-1}$ is a commutative Lie (super)algebra, whereas on the Minkowski
superspace for any $N$ and any model (except \cite{GL4}), be it
a~``conventional'', or Manin's ``exotic'' one, the tangent space at any
point possesses a natural structure of a~nilpotent Lie superalgebra.
In other words, every Minkowski superspace is a nonholonomic one,
i.e., with a nonintegrable distribution.

So we need (a) a definition of structure functions for nonholonomic
(super)spaces (this definition that solves the old Hertz-Vershik's
problem was first published in \cite{LP2}) and (b)~test which of the
coset spaces, or rather superspaces, satisfy a natural requirement: if
we throw away all odd parameters we get the conventional Einstein
equation (plus, perhaps, something else).

\newpage

In \cite{GL1} we executed this approach for every $N\leq 8$ and
several most symmetric parabolic subgroups; in our models of $N=8$
extended SUGRA, it is $Gr^8_{4}$ (dark matter?)  together with two
more copies of our Minkowski space (hell and paradise?)  that
constitute the space of extra parameters of the even, usual Minkowski
space compulsory if we wish to satisfy the above natural requirement.
These additional spaces together are analogs of ``harmonic'' space of
\cite{GI}.  Observe that the manifolds like $Gr^8_{4}$ appear in our
examples of ``distinguished'' classical spaces, the ones on which one
can write an analog of EE.

(3) The idea to apply cohomology to describe SUGRA
appeared first, perhaps, in Schwarz's paper \cite{Sz} and \cite{CDF}
but their execution of the idea is different from ours and leads
astray, we think, as far as SUGRA is concerned.

(4) Among compact Hermitian symmetric spaces, some are {\it
distinguished} by the fact that the corresponding Jordan algebra is
simple; e.g., such is the Grassmannian $Gr_{2n}^{4n}$.  In~\cite{L4}
there is given a number of examples of simple Jordan superalgebras
corresponding to simple ${\mathbb Z}$-graded Lie superalgebras of
polynomial growth.  It turns out that on manifolds locally equivalent
(in the sense of $G$-structures) to the distinguished Hermitian
spaces, one can write equations resembling the conventional EE. To
investigate how far can one stretch the analogy on supermanifolds is
an open problem.

\section{Recapitulations}

In this section we recall basic definitions~\cite{St} and retell some of
Goncharov's results~\cite{G} in a form convenient for us.

\subsection{Principal fibre bundles}

Let $M$ be a manifold of dimension $n$ over a field ${\mathbb K}$
(here: ${\mathbb R}$ or ${\mathbb C}$)
and $G$ a Lie group.
A~{\it principal fibre bundle
$P=P(M, G)$ over $M$ with group $G$} consists of a manifold $P$ and
an action of $G$ on $P$
satisfying the following conditions:

(1) $G$  acts freely on $P$ on the right;

(2) $M=P/G$ and the canonical projection $\pi :P\longrightarrow M$ is
differentiable;

(3) $P$ is locally trivial.

\begin{example}
$P=M\times G$,  the trivial bundle.  The free $G$-action on $P$ is
given by the formula $ub = (x,  ab)$ for any $u=(x,  a)\in P$,  $b\in G$.
\end{example}

\begin{example}
The bundle of linear frames over $M$.  Let $\dim M=n$.  A linear
frame $f(x)$ at a~point $x\in M$ is an ordered basis $X_1, \ldots
,  X_n$ of the tangent space $T_xM$.  Let ${\mathcal F}(M)$ be the
set of all linear frames at all points of $M$ and $\pi : {\mathcal
F}(M)\longrightarrow M$ the map such that $\pi (f(x))=x$.  The
group $GL(n)$ acts on ${\mathcal F}(M)$ on the right as follows:
if $f(x)=(X_1,  \dots ,  X_n)$ and $(a_j^i)\in GL(n)$,  then $fa=
(Y_1,  \dots ,  Y_n)$,  where $Y_i=\sum\limits_ja_i^jX_j$ is a
linear frame at $x$.  So $GL(n)$ acts freely on ${\mathcal F}(M)$
and $\pi (u)=\pi (v)$ if and only if $v=ua$ for some $a\in GL(n)$.
\end{example}

\subsection{Structure functions}
Let ${\mathcal F}(M)$ be the principal $GL(n; {\mathbb K})$-bundle
of linear frames over $M$.
Let $G\subset GL(n; {\mathbb K} )$ be a Lie group.
A $G$-{\it structure on} $M$ is a reduction
of ${\mathcal F}(M)$ to a principal $G$-bundle.

The simplest $G$-structure is the {\it flat} $G$-structure defined as
follows.  Let $V$ be ${\mathbb K} ^n$ with a fixed frame.  The flat structure
is the bundle over $V$ whose fiber over $v\in V$ consists of all
frames obtained from the fixed one under the $G$-action,  $V$ being
identified with $T_vV$.

{\bf Examples of flat structures.} The classical spaces, i.e., compact
Hermitian symmetric spaces, provide us with examples of manifolds with
nontrivial topology but flat $G$-structure.  We will shortly derive a
well-known fact that the only possible $GL(n)$-structure on any
$n$-di\-mensional manifold is always flat.

In \cite{Gu} the obstructions to identification of the $k$-th
infinitesimal neighborhood of a~point on a manifold $M$ with
$G$-struc\-ture with the $k$-th infinitesimal neighborhood of
a~point of the flat manifold $V$ with the above $G$-structure are
called {\it structure functions of order} $k$, or briefly SF.
In~\cite{Gu} and~\cite{St} it is shown that the tensors that constitute these
obstructions are well-defined provided the structure functions of
all orders $<k$ vanish.

We will write $M\sim N$ for two locally equivalent $G$-structures on
manifolds $M$ and $N$.

The classical description of the structure functions uses the notion
of the {\it Spencer cochain complex}.  It is defined as follows.  Let
$S^i$ denote the operator of the $i$-th symmetric power, prime $'$
denotes the dualization.  Set ${\mathfrak g}_{-1} = T_mM$, ${\mathfrak
g} _0 = \Lie (G)$; for $i > 0$ set:
\begin{gather*}
  {\mathfrak g}_i  = \{X\in \Hom({\mathfrak g}_{-1}, {\mathfrak g}
    _{i-1}): X(v)(w) = X(w)(v)\ \ \mbox{for any}\ \ v,
w \in {\mathfrak g} _{-1}\} \\
\phantom{{\mathfrak g}_i}{} = ({\mathfrak g}_0\otimes S^i{\mathfrak g}'_{-1})\cap
({\mathfrak g}_{-1}\otimes S^{i+1}{\mathfrak g}'_{-1}).
\end{gather*}
Finally, set $({\mathfrak g} _{-1}, {\mathfrak g} _0)_* =
\mathop{\oplus} \limits _{i \geq -1}{\mathfrak g} _i$.  It is easy to
check that $({\mathfrak g} _{-1}, {\mathfrak g} _0)_*$ possesses a
natural Lie algebra structure.  The Lie algebra $({\mathfrak g} _{-1},
{\mathfrak g} _0)_*$ is called the {\it Cartan's prolong} (the result
of {\it Cartan's prolongation}) of the pair $({\mathfrak g} _{-1},
{\mathfrak g} _0)$.

Suppose that
\begin{gather}
\mbox{the}\ \ {\mathfrak g} _0\mbox{-module}\ \ {\mathfrak g} _{-1}\ \
\mbox{is faithful}.
 \label{0.1}
\end{gather}
Then, clearly, $({\mathfrak g} _{-1}, {\mathfrak g} _0)_*\subset
{\mathfrak{vect}} (n) = {\mathfrak{der}} \, {\mathbb K} \;[x_1, \ldots
, x_n]$, where $n =\dim {\mathfrak g} _{-1}$, with
\[
{\mathfrak g}_i=\{X\in{\mathfrak{vect}}(n)_i: [X, D]\in{\mathfrak
g}_{i-1}\ \ \mbox{for any} \ \ D\in{\mathfrak g}_{-1}\}  \quad
\mbox{for}  \quad i \geq 1.
\]

Let $\Lambda^i$ be the operator of the $i$-th exterior power; set
$C^{k, s}_{({\mathfrak g} _{-1}, {\mathfrak g} _0)_*} = {\mathfrak g}
_{k-s}\otimes \Lambda^s({\mathfrak g}_{-1}')$; we often drop the
subscript of $C^{k, s}_{({\mathfrak g}_{-1}, {\mathfrak g} _0)_{*}}$
or indicate only ${\mathfrak g} _0$ since the module ${\mathfrak
g}_{-1}$ is clear.

Define the differential $\partial _s: C^{k, s}\longrightarrow
C^{k, s+1}$ by setting for any $v_1, \dots, v_{s+1}\in {\mathfrak
g}_{-1}$ (as usual, the slot with the hatted variable is to be
ignored):
\[
(\partial _sf)(v_1, \dots , v_{s+1}) = \sum_i (-1)^i \left[f\left(v_1,
\dots , \hat v _{i}, \dots , v_{s+1}\right), v_{i}\right].
\]
As expected, $\partial _s\partial _{s+1} = 0$, and the homology $H^{k,
s}_{({\mathfrak g} _{-1}, {\mathfrak g} _0)_{*}}$ of the bicomplex
$\mathop{\oplus} \limits _{k, s}C^{k, s}$ is called the $(k, s)$-th
{\it Spencer cohomology} of $({\mathfrak g} _{-1}, {\mathfrak g} _0)$.
(In the literature various gradings of the Spencer complex are in use;
ours is the most natural one.)

\begin{proposition}[\cite{Gu}]
The Spencer cohomology group
$H^{k, 2}_{({\mathfrak g}_{-1}, {\mathfrak g} _0)_{*}}$
constitutes the space of values of the structure function of order $k$.
\end{proposition}

\subsection{The case of simple
{\mathversion{bold}$({\mathfrak g} _{-1}, {\mathfrak g}_0)_*$} over
{\mathversion{bold}${\mathbb C}$}}

The following remarkable fact, though known to experts, is seldom
formulated explicitly:

\begin{proposition}\label{proposition-0.3} Let ${\mathbb K} = {\mathbb
C}$; let ${\mathfrak g} _* = ({\mathfrak g} _{-1}, {\mathfrak g}_0)_*$
be simple.  Then only the following cases are possible:

1) if ${\mathfrak g} _2 \not = 0$, then ${\mathfrak g} _*$ is either
${\mathfrak{vect}} (n)$ or its special subalgebra ${\mathfrak{svect}}
(n)$, or the subalgebra ${\mathfrak{h}} (2n)\subset {\mathfrak{vect}}
(2n)$ of hamiltonian fields;

2) if ${\mathfrak g} _2 = 0$, then ${\mathfrak g} _1 \not = 0$ and
${\mathfrak g} _*$ is the Lie algebra of the complex Lie group of
automorphisms of a compact Hermitian symmetric space.
\end{proposition}

\begin{remark} This Proposition gives a reason to impose the
restriction (\ref{0.1}) if we wish $({\mathfrak g}_{-1}, {\mathfrak
g}_0)_*$ to be simple.  On supermanifolds, where the analogue of
Proposition~\ref{proposition-0.3} does not imply similar restriction,
(or if we do not care whether ${\mathfrak g} _*$ is simple or not) we
do consider Cartan prolongs not embeddable into
${\mathfrak{vect}}(\dim {\mathfrak g} _{-1})$, see \cite{P2,P3}.
\end{remark}

Let us express Spencer cohomology in terms of Lie algebra cohomology.
Namely, observe that:
\begin{gather}
\mathop{\oplus}\limits_k H^{k, 2}_{({\mathfrak g}_{-1}, {\mathfrak
g}_0)_{*}} = H^2({\mathfrak g} _{-1}; {\mathfrak g} _*).  \label{0.2}
\end{gather}
This representation has only advantages: we loose nothing, because a
finer grading of Spencer cohomology is immediately recoverable from
the rhs of (\ref{0.2}) where it corresponds to the ${\mathbb
Z}$-grading of ${\mathfrak g} _*=({\mathfrak g}_{-1}, {\mathfrak
g}_0)_*$; moreover, there are several theorems helping to compute Lie
algebra cohomology (\cite{Fu}) whereas in order to compute Spencer
cohomology we can only use the definition.

To compute $H^2({\mathfrak g} _{-1}; {\mathfrak g} _*)$
is especially easy when ${\mathfrak g} _*$
is a simple finite dimensional Lie algebra over ${\mathbb C}$.  Indeed,
thanks to the Borel--Weil--Bott (BWB) theorem, cf.~\cite{G}, the
${\mathfrak g}_{0}$-module $H^2({\mathfrak g}_{-1}; {\mathfrak g} _*)$
has as many irreducible ${\mathfrak g}_{0}$-modules
as
$H^2({\mathfrak g}_{-1})$ which, thanks to commutativity of ${\mathfrak g} _{-1}$, is just
$\Lambda^2({\mathfrak g}_{-1}')$.  The highest weights of these irreducible modules
are also deducible from the theorem, as it is explained in~\cite{G}.
Since~\cite{G} does not give the explicit values of these weights, we
give them.  We also calculate structure functions corresponding to
case 1) of the Proposition~\ref{proposition-0.3}: we did not find these calculations
in the literature.

In what follows $R\left(\sum a_i\pi _i\right)$ denotes the irreducible ${\mathfrak g}
_0$-module (and the corresponding representation) with highest weight
$ \sum a_i\pi _i$ expressed in terms of fundamental weights as in
\cite{OV}; the weights of the ${\mathfrak{gl}} (n)$-modules,  however,  are given
for convenience with respect to the matrix units $E_{ii}$.

The classical spaces are listed in Table~1 and some of them are
baptized for convenience of further references.

Our next task is to superize Proposition~\ref{proposition-0.3} and compute the corresponding
structure functions.  For the list of ``classical'' Lie superalgebras
see~\cite{K0} (finite dimensional Lie superalgebras), \cite{GLS1}
(stringy Lie superalgebras), \cite{FLS} (Kac--Moody Lie superalgebras)
and~\cite{LSh} (or \cite{K1} and \cite{CK1,CK2}) (vectorial
Lie superalgebras). For notations of vectorial Lie superalgebras (simple and close to simple), see \cite{LSh}, \cite{GLS1}.

\begin{theorem} \label{theorem-0.1.2} 1) In case 1)~of
Proposition~\ref{proposition-0.3} the structure functions can only be
of order~$1$ {\rm (Serre, see \cite{St})}.  The actual values of
structure functions are as follows {\rm (\cite{LP2})}:

a) $H^2({\mathfrak g} _{-1}; {\mathfrak g} _*) = 0$ for ${\mathfrak g}
_* = {\mathfrak{vect}} (n)$ and ${\mathfrak{svect}} (m)$, $n, m>2$;

b) $H^2({\mathfrak g} _{-1}; {\mathfrak g} _*) = R(\pi
_3)\mathop{\oplus} R(\pi _1)$ for ${\mathfrak g} _* = {\mathfrak{h}}
(2n)$, $n>2$; $H^2({\mathfrak g} _{-1}, {\mathfrak g} _*) = R(\pi _1)$
for ${\mathfrak g} _* = {\mathfrak{h}} (4)$.

2) {\rm (Goncharov~\cite{G})} The structure functions for a space of
type $Q_3$ can be of order $3$ and constitute $R(4\pi _1)$.

The cocycles representing structure functions for a space of type $Gr_m^{n}$ (when neither
$m$ nor $n-m$ is equal to $1$, i.e., when $Gr$ is not a projective
space) belong to the direct sum of two irreducible
 (as ${\mathfrak g}_{0}$-modules) components. In this case
 ${\mathfrak g}_{0} =  {\mathfrak{sl}}(m) \oplus  {\mathfrak{sl}}(n-m)
 \oplus {\mathbb C}$ and $\Lambda^2({\mathfrak{g}}_{-1}')$ is
\begin{gather*}
\Lambda^2(({\mathbb C}^m)'\otimes ({\mathbb C}^{n-m})')=S^2({\mathbb
C}^m)'\otimes\Lambda^2 ({\mathbb C}^{n-m})'\oplus \Lambda^2({\mathbb
C}^m)'\otimes S^2({\mathbb C}^{n-m})' =\Lambda^2_+\oplus\Lambda^2_-.
\end{gather*}
The space of structure functions is the sum of two irreducible components:
the self-dual part,  $H_+$,  and antiself-dual part,  $H_-$.

The following table indicates order of components $H_{\pm}$; the
highest weight of $H_+$ (resp. $H_-$) is the sum of the highest
weights of $\Lambda^2_+$ (resp.\ $\Lambda^2_-$) and  the highest
weight of ${\mathfrak g}_{k-2}$, where $k$ is the indicated order
of the structure function:
\begin{center}
\begin{tabular}{|l|c|c|}
\hline
& SF of order $1$ & SF of order $2$\\
\hline
$m=2$, $n-m\neq 2$ & $H_-$ & $H_+$\\
$n-m=2$,  $m\neq 2$ & $H_+$ & $H_-$\\
$n-m$,  $m=2$       &       & $H_-\oplus H_+$\\
$m$,  $n-m\neq 2$   & $H_-\oplus H_+$&\\
\hline
\end{tabular}
\end{center}

The structure functions of $G$-structures of the rest of the
classical compact Hermitian symmetric spaces are the following
irreducible ${\mathfrak g}_0$-modules, where $V$ is the identity
${\mathfrak g}_0$-mo\-dule:
\begin{center}
\begin{tabular}{|c|c|c|c|}
\hline
CHSS  &$\vphantom{\Big|}  P^n$&$OGr_m$&$LGr_m$\\
\hline
conformal SF & none &$\vphantom{\Big|}
\Lambda^2\left(\Lambda^2(V')\right)\otimes V$&$
\Lambda^2\left(S^2(V')\right) \otimes V$\\
\hline
CHSS &$\vphantom{\Big|} Q_n$, $n>4$&$E_ 6^c/SO(10)\times U(1) $&$
E_7^c/E_6^c\times U(1)$\\
\hline
conformal SF & $\vphantom{\Big|} \Lambda^2V' \otimes
V$&$\Lambda^2(R(\pi _5)'))\otimes R(\pi _5)$& $ \Lambda^2(R(\pi
_1)'))\otimes R(\pi _1)$\cr \hline
\end{tabular}
\end{center}

Their order is equal to $1$ (recall that $Q_4 =
Gr_2^4$).
\end{theorem}

\subsection{Connections and structure functions}

(After \cite{M}.)  Let ${\mathcal M}$ be a supermanifold,
${\mathcal S}$ a locally free sheaf (of sections of a vector
bundle) on ${\mathcal M}$.  Locally, in a sufficiently small
neighbourhood ${\mathcal U}$, we may view ${\mathcal S}$ as a~free
module over a supercommutative superalgebra ${\mathcal F}$, which,
in the general setting, is the structure sheaf of ${\mathcal M}$.

On ${\mathcal S}$, a {\it connection} is an {\it odd} map $\nabla:
{\mathcal S} \longrightarrow {\mathcal S}\otimes_{{\mathcal F}}\Omega
^1$, where $\Omega^i = \Omega^i({\mathcal M})$ is the sheaf of
differential $i$-forms on ${\mathcal M}$.  The map $\nabla$ can be
extended to the whole de Rham complex of differential forms:
\begin{gather}
{\mathcal S} \stackrel{\nabla}{\longrightarrow } \Omega ^1
\otimes_{{\mathcal F}}{\mathcal S}\stackrel{\nabla}{\longrightarrow
}\Omega ^2\otimes_{{\mathcal F}}{\mathcal S}\cdots
\stackrel{\nabla}{\longrightarrow }\Omega ^i\otimes_{{\mathcal
F}}{\mathcal S}\cdots\label{0.4.1}
\end{gather}
by the Leibniz rule
\[
\nabla (f\otimes s) = df\otimes s + (-1)^{p(f)}f\nabla(s) \quad
\mbox{for} \quad f\in \Omega^i, \ \ s\in {\mathcal S}.
\]

Dualization determines the action of $\nabla$ on the spaces of {\it
integrable} forms, where $\Sigma _{-i}=\Hom_{{\mathcal
F}}(\Omega^i, \Vol)$ and $\Vol$ is the sheaf of volume forms:
\begin{gather}
\cdots\stackrel{\nabla}{\longrightarrow }{\mathcal
S}\otimes_{{\mathcal F}}\Sigma _{p-q-1}
\stackrel{\nabla}{\longrightarrow } {\mathcal S}\otimes_{{\mathcal
F}}\Sigma _{p-q}\stackrel{\nabla}{\longrightarrow } 0\label{0.4.1'}
\end{gather}
compatible with the $\Omega^*$-action on $\Sigma_*$ and given by the
formula
\[
\nabla(s\otimes \sigma) = T(\nabla (s))\otimes \sigma +
(-1)^{p(s)}s\otimes d(\sigma)  \quad \text{for} \quad \sigma\in
\Sigma_i, \ \ s\in {\mathcal S},
\]
where $T: \Omega^1\otimes_{{\mathcal F}} {\mathcal S}\longrightarrow
{\mathcal S}\otimes_{{\mathcal F}}\Omega^1$ is the twisting
isomorphism (mind Sign Rule).

One connection always exists: in the ${\mathcal S}$-valued de Rham
complex, set: $\nabla = d$ (more precisely, $\nabla = d\otimes
\id_{{\mathcal S}}$).  Since, as it is easy to verify, any two
connections differ by an $\Omega^*$-linear map, then any connection is of
the form $\nabla = d + \alpha$, where $\alpha\in \Omega ^1$ is called
the {\it form of the connection} or a {\it gauge field}.  We can
consider the connection as acting in the whole spaces
$\Omega^*\otimes_{{\mathcal F}}{\mathcal S}$ and ${\mathcal
S}\otimes_{{\mathcal F}}\Sigma _*$.  Then
$\nabla ^2$ $\left(=\frac{1}{2} [\nabla, \nabla]\right)$ is a
well-defined element denoted by $F_{\nabla}\in{\cal E}nd ~{\mathcal
S}\otimes_{{\mathcal F}}\Omega ^2$ and called the {\it curvature form}
of $\nabla$ or the {\it stress tensor of the Yang--Mills field}
$\alpha$.

A connection in $\Vect (M)$ is called an {\it affine} one.  An affine
connection is {\it symmetric} if
\begin{gather}
\nabla _X(Y)-  \nabla _Y(X) -[X,  Y]=0.\label{0.4.2}
\end{gather}
An affine connection $\nabla$ is called {\it compatible with the given
metric} $g$ if
\begin{gather}
g(\nabla _XY, Z)= (-1)^{p(X)p(g)}X(g(Y, Z))-(-1)^{p(X)p(Y)}g(Y, \nabla
_XZ).\label{0.4.3}
\end{gather}
Compatibility with a differential 2-form $\omega$ or another tensor of
valency $(a, b)$, say, a volume form, is similarly defined, only the
number of variable vector fields involved is not three anymore, but
$a+b+1$.

On every Riemannian manifold we have no structure functions of
order~1; hence, there always exists a unique torsion-free connection
compatible with the metric (it is called {\it the Levi--Civita
connection}) and the 2nd order structure functions are well-defined.
(This is not so for certain other $G$-structures, cf.\
Theorem~\ref{theorem-6.1} or the case of $d>1$, e.g., the case of
Minkowski {\it super}spaces.)

\subsection{Structure functions for Riemann-type structures}
In \cite{G} Goncharov considered {\it generalized conformal
structures}.  The structure functions for the corresponding
generalizations of the Riemannian structure, i.e., when Goncharov's
${\mathfrak g} _0$ is replaced with its semisimple part
$\hat{{\mathfrak g}}$ of ${\mathfrak g} = \Lie (G)$, seem to be more
difficult to compute because in these cases $({\mathfrak g} _{-1},
\hat{{\mathfrak g}} _0)_* = {\mathfrak g} _{-1}\mathop{\oplus}
\hat{{\mathfrak g}} _0$ and the BWB-theorem does not work.
Fortunately, as follows from the cohomology theory of Lie algebras, we
still have an explicit description of structure functions:

\begin{proposition}[\cite{G}, Theorem~4.7] For ${\mathfrak g} _0 =
\hat{{\mathfrak g}}$ structure functions of order $1$ are the same as
for ${\mathfrak g} _0 = {\mathfrak g}$ and structure functions of
order $2$ for ${\mathfrak g} _0 = \hat{{\mathfrak g}}$ are the same as
for ${\mathfrak g} _0 = {\mathfrak g}$ plus, additionally,
$S^2({\mathfrak g} _{-1}')$.  (Clearly, there are no structure
functions of order $>2$ for ${\mathfrak g} _0 = \hat{{\mathfrak g}}$.)
\end{proposition}

Let $G = O(n)$, i.e., $M\sim Q_{n}$.  In this case ${\mathfrak g} _1 =
{\mathfrak g} _{-1}$ and in $S^2({\mathfrak g} _{-1}')$ a
1-dimensional trivial $G$-module is distinguished; the section through
the subbundle with this subspace as a fiber is a Riemannian metric $g$
on $M$.

Let now $t$ be a structure function (the sum of its components belongs
to the distinct irreducible $O(n)$-modules that constitute
$H^2({\mathfrak g}_{-1}; {\mathfrak g} _*)$) corresponding to the
Levi--Civita connection.  The process of restoring $t$ from $g$
(compatibility condition (\ref{0.4.3})) involves differentiations thus
making any relation on $t$ into a nonlinear partial differential
equation.  Let us consider certain other restrictions on~$t$.

The values of the Riemann tensor ${\mathcal R}$ at a point of $M$
constitute an $O(n)$-module $H^2({\mathfrak g} _{-1}; {\mathfrak
g} _*)$ which contains a trivial component.  Due to complete
reducibility of finite dimensional $O(n)$-modules, we can
consider, separately, the component of ${\mathcal R}$
correspon\-ding to the trivial representation, denote it
$\mbox{Scal}$.  As is explained in
Proposition~\ref{proposition-0.3}, this trivial component is
realized as a submodule in an isomorphic copy of $S^2({\mathfrak
g} _{-1}')$, the space the metric is taken from.  Thus, we have
two matrix-valued functions, $g$ and $\mbox{Scal}$, each a~section
of the line bundle corresponding to the trivial ${\mathfrak g}
_0$-module.

What is more natural than to require their ratio to be a constant
(instead of a function)?  This condition
\begin{gather}
\mbox{Scal} = \lambda g,  \quad \mbox{where} \quad \lambda \in {\mathbb
K}, \label{EE_0}
\end{gather}
gives us ``a lesser half" of what is known as {\it Einstein
Equations} (EE).

To obtain the remaining part of EE, recall that $S^2({\mathfrak g}
_{-1}')$ consists of the two irreducible $O(n)$-components, the trivial
one and another one. A section through this other component is the
traceless Ricci tensor, $\mbox{Ric}$.  The analogs of {\it Einstein
equations} (in vacuum and with cosmological term proportional to
$\lambda$) are the {\it two} conditions: (\ref{EE_0}) and
\begin{gather}
\mbox{Ric}= 0.  \label{EE_{ric}}
\end{gather}
The remaining components of EE are invariant under conformal
transformations and do not participate in EE.

\section{Structure functions for reduced structures~--- \\
analogs of EE on manifolds}

In~\cite{G} Goncharov did not explicitly calculate the weights of
structure functions for $G$-struc\-tures corresponding to the
reduction of the generalized conformal structure.  Let us fill in
this gap: let us elucidate Proposition~\ref{proposition-0.3} for
the classical compact Hermitian symmetric spaces (CHSS).

\begin{proposition} Let ${\mathfrak g} _0$ be the semisimple part
$\hat{{\mathfrak g}}$ of ${\mathfrak g} = \Lie (G)$ corresponding to a
compact Hermitian symmetric space $X$ other than ${\mathbb O} P^2$, ${\mathbb E}$.  Then
nonconformal structure functions are all of order $2$ and as follows:
\begin{center}
\begin{tabular}{|c|c|c|}
\hline
CHSS  &$\vphantom{\Big|}{\mathbb C} P^n $&$  Gr^{m+n}_m$\\
\hline
weight of SF &$\vphantom{\Big|}R(\pi _2)$&$R(2\pi _1)\otimes R(2\pi
_1)' \mathop{\oplus} R(\pi _2)\otimes R(\pi _2)' $\\
\hline
CHSS  &$\vphantom{\Big|}OGr_m $&$  LGr _m$\\
\hline
weight of SF &$\vphantom{\Big|}R(2\pi _2)'\mathop{\oplus} R(\pi _4)'$
&$R(0, \dots , 0, -2, -2)\mathop{\oplus} R(0, \dots , 0, -4)$\\
 \hline
\end{tabular}
\end{center}
\end{proposition}

Let us show how to obtain equations similar to EE on some compact
Hermitian symmetric spaces other than $Q_n$.  Let $R$ be a section of
the vector bundle with the above structure function as the fiber; if
the space of structure function consists of two irreducible
$G$-components; denote the corresponding components of the structure
function by $R = R_1 + R_2$ in accordance with the decomposition of
the module of structure functions as indicated in the table above.  We
will consider structure functions corresponding to the canonical (in
the same sence as Levi--Civita) connection corresponding to the
$G$-structure considered.

The analogues of (\ref{EE_0})  can be defined in the following cases:

1) $Gr_{2n}^{4n}$ (turns into the conventional (\ref{EE_0}) at $n =
1$);

2) $P^{2n}$;

3) $OGr_{4n}$ (turns into the the conventional (\ref{EE_0}) at $n =
1$).

These analogues  are the equations:
\begin{gather}
v = \lambda R_2^{n} \quad (\mbox{or} \ \ v = \lambda R^{n} \ \
\mbox{if} \ \ R \ \ \mbox{has just one irreducible component}),
\label{EE_0+}
\end{gather}
where $v$ is a fixed volume element on $X$.

The analogues of (\ref{EE_{ric}}) are the equations
\begin{gather}
R_1 = 0 \quad (\mbox{if there is such a component}).
\label{EE_{ric}+}
\end{gather}
Notice that if the space of structure functions is irreducible, there
is no (\ref{EE_{ric}+}).

If structure functions of order 1 are nonzero, denote them by $T =
\mathop{\oplus} T_i$ (here the sum runs over irreducible components).
As we have quoted from~\cite{St}, the equations EE are well-defined
provided all the $T_i$ vanish.  This yields conditions similar to
Wess--Zumino constraints in SUGRA~\cite{WZ}:
\begin{gather}
T_i = 0 \quad \mbox{for every} \quad i.  \label{EE_{tor}}
\end{gather}

Notice that for all the compact Hermitian symmetric spaces the 1-st
order structure functions vanish.

Explicit computations of the structure functions for the exceptional
CHSS (see Table~1) will be given elsewhere.

\section{Analogs of EE on supermanifolds}

The theory of Lie supergroups and even Lie superalgebras is
yet new in Mathematics.  Therefore the necessary background is
gathered in a condensed form in Appendix.

We have often heard that ``the Riemannian geometry has parameters
whereas the symplectic one does not".  It is our aim to elucidate this
phrase: we have shown (Theorem~\ref{theorem-0.1.2} above,
Theorem~\ref{theorem-6.1} below and~\cite{P1}) that an {\it almost}
symplectic geometry does have parameters, the torsion, which being of
order 1 should be killed, like Wess--Zumino constraints, in order to
reduce the 2-form to a canonical form.  The curvature, alias a
structure function of order~2, might have been an obstruction to
canonical form but it vanishes.

Similar is the situation for supermanifolds.  But not quite:
${\mathfrak{o}}$ is never isomorphic to ${\mathfrak{sp}}$ whereas
the ortho-symplectic Lie superalgebra which preserves a
nondegenerate even skew-symmetric bilinear form,
${\mathfrak{osp}}^{sk}(V)$, is isomorphic to the Lie superalgebra
preserving a nondegenerate even symmetric bilinear form,
${\mathfrak{osp}} (\Pi (V))$. Still their Cartan prolongs are quite distinct:
$(V, {\mathfrak{osp}}^{sk}(V))_* = {\mathfrak{h}} (\dim (V))$ whereas
$(\Pi (V), {\mathfrak{osp}} (\Pi (V)))_* = \Pi (V)\mathop{\oplus}
{\mathfrak{osp}} (\Pi (V))$, cf.~\cite{LSh}.

Analogously, the periplectic Lie superalgebra,
${\mathfrak{pe}}^{sk}(V)$ which preserves a nondegenerate odd
skew-symmetric bilinear form is isomorphic to the Lie superalgebra,
${\mathfrak{pe}} ^{sy}(\Pi (V))$, preserving a nondegenerate odd
symmetric bilinear form; but $(V, {\mathfrak{pe}} ^{sk}(V))_* =
{\mathfrak{le}} (\dim (V))$, see \cite{LSh}, whereas $(\Pi (V),
{\mathfrak{pe}} ^{sy}(\Pi (V)))_* = \Pi (V)\mathop{\oplus}
{\mathfrak{pe}} ^{sy}(\Pi (V))$.

{\bf Possible analogues of the EE on supermanifolds with a
$G$-structure.} (Here ${\mathfrak g} = \Lie(G)$ is a simple Lie
superalgebra (${\mathbb Z}$-graded of finite growth and not
necessarily finite-dimensional) and ${\mathfrak{cg}}$ denotes the
1-dimensional trivial central extension of ${\mathfrak g}$.)

(1) The first idea is to replace ${\mathfrak{o}} (m)$ with
${\mathfrak{osp}} (m|2n)$ for a ${\mathbb Z}$-grading of the form
\[
 {\mathfrak{osp}} (m|2n) = {\mathfrak g} _{-1}\mathop{\oplus}
 {\mathfrak g} _0\mathop{\oplus} {\mathfrak g} _1 \quad
 \mbox{with} \quad {\mathfrak g} _0= {\mathfrak{cosp}} (m-2|2n) \quad
 \mbox{and} \quad m > 2.
\]

(2) The odd counterpart of this step is to replace ${\mathfrak{osp}}
(m|2n)$ with its odd (periplectic) analogues: ${\mathfrak{pe}}
^{sy}(n)$ and ${\mathfrak{spe}} ^{sy}(n)$ and the ``mixture'' of
these, ${\mathfrak{spe}} ^{sy}(n)\; +\hspace{-3.6mm} \subset {\mathbb
C} (az+bd)$, where in matrix realization we can take $d = \diag (1_n,
-1_n)$, $z = 1_{2n}$, see Appendix.

Why is $m >2$ in (1)?  If $m=2$, then ${\mathfrak g} _0=
{\mathfrak{sp}} (2n)$ and, as we know~\cite{P2}, there are no
structure functions of order 2.  Might it be that an analogue of EE is
connected not with ${\mathfrak{sp}} (2n)$, the Lie algebra of {\it
linear} symplectic transformations, but with the infinite dimensional
Lie algebra of {\it all} symplectic transformations, i.e., the Lie
algebra ${\mathfrak{h}} (2n|0)$ of Hamiltonian vector fields?
Theorem~\ref{theorem-0.1.2} states: NO. (The structure functions are
only of order~1; the corresponding eqs.  written in~\cite{P1}, though
interesting, do not resemble EE.)

Let us not give up: the Lie algebra ${\mathfrak{o}} (m)$, as well as
${\mathfrak{h}} (2n|0)$, has one more analogue~--- the Lie
superalgebra ${\mathfrak{h}} (0|m)$ of Hamiltonian vector fields on
$(0|m)$-dimensional super\-ma\-ni\-fold.  So other possibilities are:

(3) replace ${\mathfrak{osp}} (m|2n)$ with ${\mathfrak{h}} (2n|m)$,
where $m \not = 0$.

Since we went that far,  let us go further still and

(4) replace ${\mathfrak{h}} (2n|m)$ in (3) with ${\mathfrak{k}}
(2n+1|m)$; and, moreover, consider ``odd'' analogues of (3) and (4):

(5) replace ${\mathfrak{pe}} (n)$ and ${\mathfrak{spe}} (n)$ in (2)
with ${\mathfrak{le}} (n)$ and ${\mathfrak{sle}} (n)$, ${\mathfrak m}
(n)$ or ${\mathfrak b}_{\lambda}(n)$.  For the definition of these and
other simple vectorial Lie superalgebras see \cite{Shch}.

In the next sections we will list structure functions for some of
the possibilities (1)--(6). The remaining ones are {\bf open problems}.

\begin{remark} One should also investigate the cases associated with
${\mathbb Z}$-grading of Kac--Moody (twisted loop) superalgebras of
the form $\mathop{\oplus}\limits _{|i|\leq 1} {\mathfrak g} _i$.
Nobody explored yet this infinite dimensional possibility.  Clearly,
there are ``trivial'' analogues of compact Hermitian symmetric spaces,
namely, the manifolds of loops with values in any finite dimensional
compact Hermitian symmetric space.  Remarkably, there are also
``twisted'' versions of these compact Hermitian symmetric spaces
associated with twisted loop algebras and superalgebras,
cf.~\cite{LSV1}.  It is not known, however, how to calculate the
cohomology of Kac--Moody algebras with this type of coefficients, even
in the ``trivial cases''.
\end{remark}

\section{Spencer cohomology of {\mathversion{bold}${\mathfrak{osp}}
(m|n)$}}

{\bf {\mathversion{bold}${\mathbb Z}$}-gradings of depth 1.} All these
gradings are of the form ${\mathfrak g} _{-1}\mathop{\oplus}
{\mathfrak g} _0 \mathop{\oplus} {\mathfrak g} _1$ with ${\mathfrak g}
_1 \simeq {\mathfrak g} _{-1}'$ as ${\mathfrak g} _0$-modules.

\begin{proposition}[\cite{K, LSV1}]
For ${\mathbb Z}$-gradings of depth~$1$ of ${\mathfrak{osp}} (m|2n)$
the following cases are possible:

a) ${\mathfrak{cosp}} (m-2|2n)$, ${\mathfrak g} _{-1} = \id$;

b) ${\mathfrak{gl}}(r|n)$ if $m = 2r$, ${\mathfrak g} _{-1} =
 \Lambda^2(\id)$.
\end{proposition}

{\bf Cartan prolongs of {\mathversion{bold}$({\mathfrak g} _{-1},
{\mathfrak g} _0)$} and {\mathversion{bold}$({\mathfrak g} _{-1},
\hat{{\mathfrak g} }_0)$}.}

\begin{proposition} a) $({\mathfrak g} _{-1}, {\mathfrak g}
_0)_*={\mathfrak g}$ except for the case Proposition~5b)  for $r =
3$, $n = 0$ when $({\mathfrak g}_{-1}, {\mathfrak g}_0)_*=
{\mathfrak{vect}} (3|0)$.

b) $({\mathfrak g}_{-1}, \hat{{\mathfrak g}}_0)_*= {\mathfrak g}
_{-1}\mathop{\oplus}\hat{{\mathfrak g}}_0$.
\end{proposition}

{\bf Structure functions.} Cases a) and b) below correspond to
cases of ${\mathbb Z}$-gradings from Proposition~5.  The cases $mn
= 0$ are dealt with in~\cite{G} and Introduction.

\begin{theorem} \label{thorem-3.3} a) As $\hat{{\mathfrak
g}}_0$-module, $H^{2, 2}_{({\mathfrak g}_{-1}, \hat{{\mathfrak g}
}_0)_{*}} = S^2(\Lambda^2({\mathfrak g}_{-1}))/\Lambda^4({\mathfrak g}
_{-1})$ and splits into the direct sum of three irreducible components
whose weights are given in Table~3.

As ${\mathfrak g} _0$-module, $H^{2, 2}_{({\mathfrak g}_{-1},
{\mathfrak g} _0)_{*}} = H^{2, 2}_{({\mathfrak g}_{-1}, \hat{{\mathfrak
g}}_0)_{*}}/S^2({\mathfrak g}_{-1})$. It is irreducible and
its highest weight is given in Table~3.
For $k\neq 2$ structure function
vanish.

b) If $r\neq n$, $n+2$, $n+3$, then
$H^2({\mathfrak g}_{-1}; {\mathfrak g}_*)$
is an irreducible ${\mathfrak g}_0$-module
and its highest weight is given in Table~4.

The cases $r = 4$, $n = 0$ and $r = 2$, $n = 1$ coincide, respectively,
with the cases considered in a) for ${\mathfrak{o}} (8)$ and
${\mathfrak{osp}} (4|2)$.
\end{theorem}

\section{Spencer cohomology of {\mathversion{bold}${\mathfrak{spe}}
(n)$}}

\begin{proposition}[Cf.~\cite{K} with \cite{LSV1}] All ${\mathbb
Z}$-gradings of depth $1$ of ${\mathfrak g}$ are listed in Table~1
of~\cite{LSV1}.  They are:

a)${\mathfrak g} _0 = {\mathfrak{sl}} (m|n-m)$,
${\mathfrak g} _{-1}= \Pi(S^2(\id))$, ${\mathfrak g} _1=
\Pi(\Lambda^2(\id'))$;

b) ${\mathfrak g} _0 = \langle \tau + (n - 1)z\rangle
 \; +\hspace{-3.6mm}\supset {\mathfrak{spe}}(n-1),
{\mathfrak g} _{-1} = \id, {\mathfrak g} _{1} = \id' = \Pi (\id),
{\mathfrak g} _{1} = \Pi (\langle {\bf 1}\rangle)$.

Here $\tau = \diag (1_{n-1}, -1_{n-1})$,
$z = 1_{2n-2}$,
the sign ${\mathfrak a} \; +\hspace{-3.6mm}\supset {\mathfrak b}$ denotes
a semidirect sum of Lie superalgebras, the ideal is on the right,
$\id$ is  endowed with a nondegenerate supersymmetric odd bilinear form.
In these cases ${\mathfrak g} _* = {\mathfrak g}$.

If ${\mathfrak g} _0 = {\mathfrak{cpe}}(n-1),
{\mathfrak g} _{-1} = \id,$ then ${\mathfrak g} _* = {\mathfrak{pe}}(n)$.
If ${\mathfrak g} _0 = {\mathfrak{spe}}(n-1),
{\mathfrak{pe}}(n-1)$ or ${\mathfrak{cspe}}(n-1)$ and
${\mathfrak g} _{-1} = \id$, then
${\mathfrak g} _* = {\mathfrak g} _{-1}\oplus {\mathfrak g} _0$.

\end{proposition}

\begin{theorem} \label{theorem-4} a) The nonvanishing structure
function are of order $1$, and in the cases when they
constitute a completely reducible
${\mathfrak g} _0$-module, the corresponding
highest weights are given in Table~5.

b) For ${\mathfrak g} _0 = {\mathfrak{spe}} (n-1)$, ${\mathfrak{pe}} (n-1)$, ${\mathfrak{cspe}} (n-1)$,$\langle \tau + (n - 1)z\rangle
\; +\hspace{-3.6mm}\supset {\mathfrak{spe}}(n-1)$,
${\mathfrak{cpe}} (n-1)$, and ${\mathfrak g} _{-1} = \id$
all structure functions vanish except for $H^{1, 2}_{{\mathfrak{spe}} (n-1)} = \Pi (\id) =\Pi (V_{\varepsilon _{1}})$ and there are the following nonsplit
exact sequence of ${\mathfrak{spe}} (n-1)$-modules:
(here $\varepsilon _{1}, \ldots, \varepsilon _{n-1}$ is the standard
basis of the dual space  to the space of diagonal matrices in
${\mathfrak{pe}} (n-1)$. $V_{\lambda}$~denotes the irreducible
${\mathfrak{pe}} (n-1)$-module with highest weight $\lambda$
and even highest vector)
\begin{gather*}
0\longrightarrow V_{\varepsilon _{1}+\varepsilon _{2}} \longrightarrow
H^{2, 2}_{{\mathfrak{spe}} (n-1)} \longrightarrow \Pi (V_{2\varepsilon
_{1}+ 2\varepsilon_{2}}) \longrightarrow 0 \quad \mbox{for} \quad n>4,\\
0 \longrightarrow X\longrightarrow H^{2, 2}_{{\mathfrak{spe}} (3)}
\longrightarrow \Pi (V_{3\varepsilon _{1}})\longrightarrow 0,
\end{gather*}
where $X$ is determined from the
following nonsplit exact sequences of ${\mathfrak{spe}} (3)$ -modules:
\[
0 \longrightarrow V_{\varepsilon _{1}+\varepsilon _{2}}
\longrightarrow X \longrightarrow \Pi (V_{2\varepsilon _{1}+2\varepsilon
_{2}}) \longrightarrow 0.
\]
Also
\[
0 \longrightarrow H^{2, 2}_{{\mathfrak{spe}(n-1)}}
 \longrightarrow H^{2, 2}_{{\mathfrak{pe}(n-1)}}
\longrightarrow
V_{2\varepsilon _{1}} \longrightarrow 0,  \quad if  \quad n > 3;
\]
\[
0 \longrightarrow H^{2, 2}_{{\mathfrak{spe}(n-1)}}
 \longrightarrow H^{2, 2}_{{\mathfrak{cspe}(n-1)}}
\longrightarrow
V_{2\varepsilon _{1}} \longrightarrow 0,  \quad if  \quad n > 3;
\]
$H^{2, 2}_{{\langle \tau + (n - 1)z\rangle
 \; +\hspace{-2mm}\supset\; {\mathfrak{spe}}(n-1)}} = \Pi
(V_{2\varepsilon _{1}+2\varepsilon _{2}})$ is an irreducible
 ${\mathfrak{spe}(n-1)}$-module if $n>4$
 and
\[
0 \longrightarrow \Pi (V_{2\varepsilon _{1}+2\varepsilon _{2}})
\longrightarrow H^{2, 2}_{{\langle \tau + 3z\rangle \;
+\hspace{-2mm}\supset\; {\mathfrak{spe}}(3)}} \longrightarrow
 \Pi (V_{3\varepsilon _{1}})\longrightarrow 0.
\]
Finally,
\[
0 \longrightarrow H^{2, 2}_{{\langle \tau + (n-1)z\rangle \;
+\hspace{-2mm}\supset \; {\mathfrak{spe}}(n-1)}} \longrightarrow
H^{2, 2}_{{\mathfrak{cpe}(n-1)}}\longrightarrow V_{2\varepsilon
_{1}}\longrightarrow 0  \quad if \quad n > 3.
\]
Moreover, $ H^{2, 2}_{{\mathfrak{cpe}(n-1)}} =
\Pi(S^2(\Lambda^2(\id)/ \Pi(\langle 1\rangle))/\Lambda^4(\id)).$

\end{theorem}

\section{An analogue of a theorem by Serre:\\
on involutivity of {\mathversion{bold}${\mathbb Z}$}-graded Lie
superalgebras}

{\samepage Theorem~\ref{theorem-0.1.2}, part of which we have
attributed above to Serre, is actually a corollary of Serre's
initial statement~\cite{St}. Before we formulate it, recall that
the notion of involutivity comes from very practical problems: how
to solve differential equation with the help of
a~computer~\cite{LS}. Let $\pi: E\longrightarrow B$ be the bundle
(it sufficies to consider  the trivial bundle with base
$B={\mathbb R}^n$ and the fiber ${\mathbb R}^m$); let $J_qE$ be
the space of $q$-jets of sections of the bundle $\pi$. Let
$V(E)\subset TE$ be the vertical bundle, i.e., the kernel of the
map $T\pi$. With every system of differential equations
$DE_q\subset J_qE$ of order $q$ in $m$ unknown functions of $n$
variables we can associate a subbundle $N_q\subset V^{(q)}J_qE$,
where $V^{(q)}J_qE$ is the vertical bundle with respect to the
projection $\pi^q_{q-1}: J_qE\longrightarrow J_{q-1}E$ as a
subbundle of $S^q(T^*B)\otimes VE$. The subbundle $N_q$ is called
the {\it geometric symbol} of the system $DE_q$. Set
\[
N_q^{(s)}:=\left\{ f\in N_q\Big| \frac{\partial f}{\partial x_i}=0
 \text{ for $i=1, \dots, s$}\right\},
\]
where $f\in{\mathbb R}^m$ and the derivatives are taken coordinate-wise.
Let $P^m_q$ be the $m$th tensor (symmetric, actually)
 power of the space of degree $\leq q$ polynomials (in $n$ variables).
  The first prolongation of $N_q$ is defined to be
\[
N_{q+1}:=\left\{ f\in P^m_{q+1}\Big| \frac{\partial f}{\partial
x_i}\in N_q \text{ for $i=1, \dots, n$}\right\}. \] The symbol
$N_q$ is said to be {\it involutive} if
\[
\dim N_{q+1}=\dim N_q +\dim N_q^{(1)}+\cdots +\dim N_q^{(n-1)}.
\]
(usually, the lhs is smaller).}

Similarly, let $g\subset {\rm Hom}(V, W)$ be a subspace and $g^{(i)}$ the $i$th Cartan prolongation of $g$ (defined above for $W=V$). For any subspace $H\subset V$ set:
\[
g_{H}:=\{F\in g\mid F(h)=0\text{ for any }h\in H\}.
\]
Let $r_k=\mathop{\min}\limits_{\dim H=k}g_{H}$.
It is not difficult to show that
\begin{gather}
\dim g^{(1)}\leq r_0+r_1+\dots +r_{k-1}. \tag{$*$}
\end{gather}
The space $g$ is called {\it involutive} if there is an equality
in $(*)$. It is not difficult to see that if $g$ is involutive,
then $g^{(1)}$ is also involutive. So, speaking about Lie algebras
which are Catran prolongs it suffices to consider involutivity of
their linear parts.

Let ${\mathfrak g} = \mathop{\oplus} \limits _{ k \geq -1} {\mathfrak
g} _k$ be a ${\mathbb Z}$-graded Lie algebra, $\{a_1, \ldots,
a_n\}$ be a basis of ${\mathfrak g} _{-1}$.  Clearly, the map
\[
\ad_{a_r}: {\mathfrak g} \longrightarrow  {\mathfrak g} ,
  \quad x \mapsto [a_r,  x]
\]
is a homomorphism of ${\mathfrak g} _{-1}$-modules.
  In accordance with the above, we say that a ${\mathbb
Z}$-graded Lie algebra of the form ${\mathfrak g}
=\mathop{\oplus}\limits_{ k \geq -1} {\mathfrak g} _k$ is called {\it
involutive} if all the maps $\ad_{a_r}$ are onto.

Serre observed that involutivity property considerably simplifies
computation of cohomology: if ${\mathfrak{g}}_*$ is involutive, then for every $i$
\begin{gather}
H^i({\mathfrak{g}}_-; {\mathfrak{g}}_*) \text{ is supported in the
lowest possible degree.}\tag{$**$}
\end{gather}

To superize the notion of involutivity, we have to require surjectivity
 of the maps $\ad_{a_r}$ for $a_r$ even.
  Additionally we must demand
vanishing of the homology with respect to each differential
 given by the odd map $\ad_{a_r}$ (the homology is
well-defined thanks to the Jacobi identity).
More precisely, for any Lie
superalgebra ${\mathfrak g} =\mathop{\oplus}\limits_{ k \geq -1}
{\mathfrak g} _k$ set:
\[
{\mathfrak g} ^{r } = \ker \, \ad_{a_1}\cap \ker \, \ad_{a_2}\cap
\cdots \cap \ker \,\ad_{a_{r }}.
\]
Clearly, ${\mathfrak g} ^{r } = \mathop{\oplus} \limits_{k \geq
-1}{\mathfrak g} ^{r }_k$, where ${\mathfrak g} ^{r }_k={\mathfrak g}
^{r }\cap {\mathfrak g}_k$.  Notice that $\ad_{a_r}\left({\mathfrak g}
^{r-1}_k\right) \subset {\mathfrak g} ^{r-1}_{k-1}$.  The Lie
superalgebra ${\mathfrak g} =\mathop{\oplus}\limits_{ k \geq -1}
{\mathfrak g} _k$ will be called {\it involutive} if the following
conditions are fulfilled:

(1) ${\mathfrak g} ^n  = {\mathfrak g} _{-1}$ (recall that $n=\dim {\mathfrak
g} _{-1}$);

(2) $\ad_{a_r}\left({\mathfrak g} ^{r-1}\right) = {\mathfrak g}
^{r-1}$ if $a_r$ is even;

(3) $\ad_{a_r}\left({\mathfrak g} ^{r-1}\right) = {\mathfrak g} ^{r }$
if $a_r$ is odd.

The cohomology group $H^{i}({\mathfrak g}_{-1};  {\mathfrak g})$
has a natural ${\mathbb Z}$-grading:

$H^{i}({\mathfrak g}_{-1};  {\mathfrak g}) = \mathop{\oplus}_{k \geq -1}
H^{i, k}({\mathfrak g}_{-1};  {\mathfrak g})$
 induced by the
${\mathbb Z}$-grading of ${\mathfrak g}$.

\begin{theorem} \label{theorem-5} {\em (\cite{P1})}
Let ${\mathfrak g}$ be involutive.
Then
if $i \geq 0$ and $k \geq 0$, then
 $H^{i, k}({\mathfrak g}_{-1};  {\mathfrak g}) = 0$.
\end{theorem}

\section{Spencer cohomology of vectorial Lie superalgebras\\
in their standard grading}

\begin{theorem}[cf. Theorem~\ref{theorem-0.1.2}]\label{theorem-6.1}

1) For ${\mathfrak g} _* = {\mathfrak{vect}} (m|n)$ and
${\mathfrak{svect}} (m|n)$ the structure functions vanish except for
${\mathfrak{svect}} (0|n)$ when the structure functions are of order
$n$ and constitute the ${\mathfrak g}_0$-module $\Pi^n({\bf 1})$.

2) For ${\mathfrak g} _* = {\mathfrak{h}} (0|m)$ for $m >4$, and
${\mathfrak g} _* = {\mathfrak{h}} (2n|m)$ for $m n\neq 0$, the
nonzero structure functions are $\Pi (R(3\pi _1)\mathop{\oplus}
R(\pi _1))$ of order~$1$.

3) For ${\mathfrak g} _* =  {\mathfrak{h}}^o (0|m)$, $m >4$, the nonzero
structure functions are same as for ${\mathfrak{h}} (0|m)$ plus an
additional direct summand $\Pi ^{m-1}(R(\pi _1))$ of order $m-1$.

4) For ${\mathfrak g} _* = {\mathfrak{sle}}(n)$, $n >1$, the nonzero
structure functions are $H^{1, 2}_{{\mathfrak{sle}} (n)} =
S^3({\mathfrak g} _{-1})$, $H^{2, 2}_{{\mathfrak{sle}} (n)} = \Pi
({\bf 1})$, $H^{n, 2}_{{\mathfrak{sle}} (n)} = \Pi ^n({\bf 1})$.
\end{theorem}

Thus, on almost symplectic manifolds with nondegenerate and non-closed
form $\omega$, there is an analog of torsion --- structure function of
order~1, namely $d\omega$.  Since the space of 3-forms splits into the
space of forms proportional to $\omega$ and the complementary space of
``primitive'' forms, there are two components of this torsion:
$d\omega=\lambda \omega+P$.  If the primitive component vanishes, we
have a nice-looking equation:
\begin{gather}
d\omega=\lambda \omega \quad\mbox{for some} \quad
\lambda\in\Omega^1.\label{6.1.1}
\end{gather}
The other component of ``torsion'' must also vanish for the supermanifold
to be symplectic, not almost symplectic.

{\bf An analogue of Einstein equation on almost periplectic
supermanifolds.} Let~$\omega _1$ be the canonical odd 2-form and $R$ a
$2$-form which is a section through $H^{2, 2}_{{\mathfrak{sle}} (n)} =
\Pi ({\bf 1})$.  This gives rise to an analogue of (\ref{EE_0}) for
${\mathfrak{sle}} (n)$:
\begin{gather}
\omega _1 = \lambda R,  \quad   \lambda \in {\mathbb C}.
\label{EE_0-sle}
\end{gather}

The equation (\ref{EE_0-sle}) are well-defined {\it provided the
irreducible components of the $1$st order structure functions, the
elements from $H^{1, 2}_{{\mathfrak{sle}} (n)}$ vanish}.  Denote by
${\bf Tor}_i$ ($i=1, 2$) the components of the torsion tensor; then
these conditions are:
\[
{\bf Tor} _i=0  \quad \mbox{for}  \quad i=1, 2.
\]
When the torsion components vanish, we can reduce the nondegenerate
odd 2-form to the canonical form (cf.~\cite{L0} and \cite{Sh}).

\section{Proof of Theorem~\ref{theorem-5}}

{\bf The long exact sequence.} Let ${\mathfrak g}$ be a Lie
superalgebra and
\begin{gather}
0\longrightarrow A\overset{\partial _0}{\longrightarrow}
C\overset{\partial _1}{\longrightarrow}
B\longrightarrow 0, \nonumber\\
 \mbox{where} \quad p(\partial _0)=\bar{0}
 \quad \mbox{and} \quad \partial _1  \quad \mbox{is either even or
odd}, \label{*}
\end{gather}
be a short exact sequence of ${\mathfrak g}$-modules.  Let $d$ be the
differential in the standard cochain complex of the Lie superalgebra
${\mathfrak g}$, cf.~\cite{Fu}.

Consider the long  sequence of cohomology:
\begin{gather}
\cdots\overset{\partial }{\longrightarrow} H^i({\mathfrak g} ;
A)\overset{\partial _0}{\longrightarrow} H^i({\mathfrak g} ;
C)\overset{\partial _1}{\longrightarrow} H^i({\mathfrak g} ;
B)\overset{\partial }{\longrightarrow} H^{i+1}({\mathfrak g} ; A)
\overset{\partial _0}{\longrightarrow} \cdots, \label{**}
\end{gather}
where $\partial _i $ is the differential induced by the namesake
differential in (\ref{*}), and $\partial =\partial _0^{-1}\circ d\circ \partial _1^{-1}$.
Since $ \partial _0$ and $\partial _1 $ commute with $d$, the sequence
(\ref{**}) is well-defined and the same arguments as for Lie
algebras~\cite{Fu} demonstrate that the long sequence of cohomology
(\ref{**}) induced by (\ref{*}) is exact.

Thus, the short exact sequences of ${\mathfrak g} _{-1}$-modules,
where $\partial _1 =\ad_{a_r}$ and $ \partial _0$ is the embedding
${\mathfrak g} ^{r }\subset {\mathfrak g} ^{r-1}$:
\begin{gather*}
0\longrightarrow {\mathfrak g} ^{r}\overset{\partial _0}{\longrightarrow}
{\mathfrak g} ^{r-1}\overset{\partial _1}{\longrightarrow}
{\mathfrak g} ^{r-1}\longrightarrow 0 \quad \mbox{for} \quad p(a_r)=\bar{0},\\
0\longrightarrow {\mathfrak g} ^{r}\overset{\partial _0}{\longrightarrow}
{\mathfrak g} ^{r-1}\overset{\partial _1}{\longrightarrow}
{\mathfrak g} ^{r }\longrightarrow 0 \quad \mbox{for} \quad p(a_r)=\bar{1}
\end{gather*}
induce the long exact sequences of
cohomology
\begin{gather}
\cdots\longrightarrow H^i\left({\mathfrak g} _{-1}; {\mathfrak g} ^{r}\right)
\overset{\partial _0}{\longrightarrow}
H^i\left({\mathfrak g} _{-1}; {\mathfrak g} ^{r-1}\right)
\overset{\partial _1}{\longrightarrow} H^i\left({\mathfrak g}
_{-1}; {\mathfrak g} ^{r-1 }\right)\nonumber\\
\overset{\partial }{\longrightarrow} H^{i+1}\left({\mathfrak g} _{-1};
{\mathfrak g} ^{r}\right) \longrightarrow \cdots, \label{a}\\
\cdots\longrightarrow H^i\left({\mathfrak g} _{-1}; {\mathfrak g} ^{r}\right)
\overset{\partial _0}{\longrightarrow}
H^i\left({\mathfrak g} _{-1}; {\mathfrak g} ^{r-1}\right)
\overset{\partial _1}{\longrightarrow} H^i\left({\mathfrak g}
_{-1}; {\mathfrak g} ^{r}\right)\nonumber\\
 \overset{\partial }{\longrightarrow}
 H^{i+1}\left({\mathfrak g} _{-1}; {\mathfrak g}
^{r}\right) \longrightarrow \cdots. \label{b}
\end{gather}

\begin{lemma} In   sequences (\ref{a}) and  (\ref{b}) we have $\partial _1=0$.
\end{lemma}

\begin{proof} First, consider sequence  (\ref{a}). Let $f\in C^i\left({\mathfrak g} _{-1}; {\mathfrak g}
^{r-1}\right) $ and $df=0$.  Define $c\in C^{i-1}\left({\mathfrak g}
_{-1}; {\mathfrak g} ^{r-1}\right)$ by the formula
\[
c(x_1,  \dots ,  x_{i-1})=f(x_1,  \dots ,  x_{i-1},  a_r).
\]
Then (here the sign $\pm$ is determined by Sign Rule)
\begin{gather*}
dc(x_1,  \dots ,  x_{i })=\sum \pm \ad_{x_j}c(x_1,  \dots , \hat x_j,
\dots ,
 x_{i })\\
 {}=\sum \pm \ad_{x_j}f(x_1,  \dots ,  \hat x_j,  \dots ,
 x_{i },   a_r)=
df(x_1,  \dots ,  x_{i },  a_r)-\pm \ad_{a_r}f(x_1,  \dots ,  x_{i }).
\end{gather*}
Since $df=0$,  we have $dc= \pm a_rf$.  Thus,  $a_rf$ is exact,
and,  therefore,  $\partial _1=0$.

Let us prove now  (\ref{b}). Let $l=\mathbb C a_r$.
Then $H^p(l; {\mathfrak g} ^{r
-1}) =0$ for $p>0$ and $H^0(l; {\mathfrak g} ^{r
-1}) = {\mathfrak g} ^{r}$. Since ${\mathfrak g} _{-1}=
l\oplus {\mathfrak g} _{-1}/l$, the spectral sequence with respect
 to the ideal $l\subset {\mathfrak g} _{-1}$ immediately gives
\[
H^{i}({\mathfrak g} _{-1};  {\mathfrak g} ^{r })\simeq H^{i}({\mathfrak g} _{-1}/l;  {\mathfrak g} ^{r })
\]
Therefore, $H^{i}({\mathfrak g} _{-1};  {\mathfrak g} ^{r-1})\subset H^{i}({\mathfrak g} _{-1};  {\mathfrak g} ^{r })$ and $\partial_0$ is surjective. Thus, $\partial_1=0$. \end{proof}

\begin{corollary} The long exact sequences (\ref{a}) and (\ref{b})
can be reduced to the following short exact sequences
\begin{gather}
0\longrightarrow H^{i-1}\left({\mathfrak g} _{-1}; {\mathfrak g}
^{r-1}\right) \overset{\partial}{\longrightarrow}
H^i\left({\mathfrak g} _{-1}; {\mathfrak g} ^{r
}\right)\overset{\partial _0}{\longrightarrow} H^i\left({\mathfrak g}
_{-1}; {\mathfrak g} ^{r-1 }\right)\overset{\partial _1}{\longrightarrow}
0,\label{a'}\\
0\longrightarrow H^{i-1}\left({\mathfrak g} _{-1}; {\mathfrak g}
^{r}\right) \overset{\partial}{\longrightarrow} H^i\left({\mathfrak
g} _{-1}; {\mathfrak g} ^{r }\right)\overset
{\partial_0}{\longrightarrow} H^i\left({\mathfrak g} _{-1}; {\mathfrak g}
^{r-1}\right) \overset{\partial _1}{\longrightarrow}0.  \label{b'}
\end{gather}
\end{corollary}

Now we can prove the theorem by induction on $r$.  First of all,
${\mathfrak g} ^n={\mathfrak g} _{-1}$ by condition~(1) of
involutivity.  So we have ${\mathfrak g} ^n_k=0$ for $k \geq 0$ and
\[
H^{i, k}\left({\mathfrak g} _{-1}; {\mathfrak g} ^{n}\right)=0  \quad
\mbox{if} \quad i \geq 0 \quad \mbox{and} \quad k \geq 0.
\]

Then consider the term of degree $k$ in (\ref{a'}) and (\ref{b'}).  We
obtain the exact sequences
\begin{gather}
0\longrightarrow H^{i-1, k}\left({\mathfrak g} _{-1}; {\mathfrak g}
^{r-1}\right) \overset{\partial}{\longrightarrow} H^{i,
k}\left({\mathfrak g} _{-1}; {\mathfrak g} ^{r
}\right)\overset{\partial _0}{\longrightarrow} H^{i,
k}\left({\mathfrak g} _{-1}; {\mathfrak g} ^{r-1 }\right)
\overset{\partial _1 }{\longrightarrow}0,\label{a''}\\
0\longrightarrow H^{i-1, k}\left({\mathfrak g} _{-1}; {\mathfrak g}
^{r}\right) \overset{\partial}{\longrightarrow} H^{i,
k}\left({\mathfrak g} _{-1}; {\mathfrak g} ^{r
}\right)\overset{\partial _0}{\longrightarrow} H^{i,
k}\left({\mathfrak g} _{-1}; {\mathfrak g} ^{r-1}\right)
\overset{\partial _1 }{\longrightarrow} 0.  \label{b''}
\end{gather}
It follows immediately from (\ref{a''})
 for $p(a_r)=\bar{0}$ and from
(\ref{b''}) for $p(a_r)=\bar{1}$
that $H^{i, k}\left({\mathfrak g} _{-1}; {\mathfrak g}^{r-1 }\right)=0$.
 The theorem is proved.

\section{Open problems: Riemann tensors\\ on curved
supergrassmannians}

Denote by ${\mathfrak g} (m|n)$ either of the Lie superalgebras
${\mathfrak{h}} (2m|n)$, $ {\mathfrak{h}}^o (n)$ or ${\mathfrak{k}}
(2m+1|n)$; let $F(m|n-2)$ be the superspace of ``functions'' which in
our case are polynomials or power series on which ${\mathfrak g} _0$
naturally acts.

In \cite{LSV1}, Table 5, there are listed all ${\mathbb Z}$-gradings of
${\mathfrak g}={\mathfrak g} (m|n)$ of the form ${\mathfrak g} =
{\mathfrak g} _{-1} \mathop{\oplus} {\mathfrak g} _{0}\mathop{\oplus}
{\mathfrak g} _1$.  For them, ${\mathfrak g} _1 \simeq{\mathfrak g}
_{-1}'$, ${\mathfrak g} _{-1} = F(m|n)$, ${\mathfrak g}_{0}
={\mathfrak g} (m|n-2)\mathop{\oplus} F(m|n-2)$ for $n>1$, and if
$n>2$, then ${\mathfrak g} _{-1 }$ is not purely odd and is isomorphic
to the tangent space to the total space of the what is called Fock
bundle over a $(2m|n-2)$-dimensional symplectic supermanifold or its
version for the contact supermanifold.

In 1985 Yu~Kochetkov informed us that he showed (unpublished) that for
${\mathfrak g} (m|n) = {\mathfrak{h}} (2m|n)$ or ${\mathfrak{h}}
(0|n)$ there is always a trivial component (perhaps, there are
several) in the space of 2nd order structure functions for the Riemann-like tensors, so there are
analogues of (\ref{EE_0}). Observe, that for Weyl-like tensors (for conformal structures) there is no trivial modules and, this is expected since trivial modules correspond to filtered deformations, cf. e.g., \cite{CK1}.

One of us (EP) managed to calculate structure functions of order 1 for
${\mathfrak g} =  {\mathfrak{h}}^o (0|6)$.  The space of these structure
functions is nonzero; in addition, it is not completely reducible,
some of the indecomposable components look as complicated as follows,
where $x$ and $y$ are some irreducible components (the same symbol
denotes an isomorphic copy):
\begin{center}
\begin{tabular}{ccccc}
$x $&$\longrightarrow  $&$y $&$\longleftarrow  $&$  x $\\
$\downarrow$& &$\uparrow$&& $\downarrow$\\
$y $&$\longleftarrow   $&$   x$&$\longrightarrow  $&$y $\\
$\uparrow$& &$\downarrow$&&$\uparrow$\\
$   x $&$\longrightarrow  $&$y $&$\longleftarrow  $&$   x $
\end{tabular}
\end{center}

Since structure functions of order 1 must vanish in order for the
analogues of EE be well-defined, these structure functions constitute
constraints similar to the Wess--Zumino constraints in supergravity.
Here we encounter an amazing situation: the lack of complete
reducibility implies that only part of these constraints (depicted by
$x$) are relevant.

We have no idea how to approach analytically other, especially
infinite dimensional, cases: the number of structure function grows
quickly with $m$ and $n$.  The only way we see at the moment is to arm
ourselves with computers, e.g., Grozman's package SuperLie~\cite{GL2}.
A first result in this direction is calculation of structure functions
for curved supergrassmannian ${\mathcal G}_{0|2}^{0|4}$ and its
``relatives'' resulting in an unconventional and unexpected version of
supergravity equations~\cite{GL4}.

\newpage

\section{Tables}

{\bf Notations in tables.}  We use the notational conventions
of~\cite{S} and definitions adopted there.

In Table 1: ${\mathfrak{s}} = (\Lie(S_c))\otimes {\mathbb C}$,
NCHSS is an abbreviation for {\it noncompact Hermitian symmetric
space}, in the diagram of ${\mathfrak{s}}$ the maximal parabolic
subalgebra ${\mathfrak p} = \Lie(P)$, such that $X$ can be
represented as $(S_c)^{{\mathbb C}}/P$, is determined by one
vertex: the last one in cases 0, 2, 3, $E_7$, $E_8$, the first one
in case 4, the $p$-th one in case 1.  Symbol ${\mathfrak{cg}}$
denotes the trivial central extension of the Lie (super) algebra
${\mathfrak g}$.

In Table 2: we call a homogeneous space $G/P$, where $G$ is a
simple Lie supergroup~$P$ its parabolic subsupergroup
corresponding to several omitted generators of a Borel subalgebra
(description of these generators can be found in~\cite{GL3}), of
{\it depth} $d$ and {\it length} $l$ if such are the depth and
length of $\text{Lie} (G)$ in the ${\mathbb Z}$-grading compatible
with that of $\text{Lie}(P)$. Note that all superspaces of Table~2
possess an hermitian structure (hence are of depth~1) except
$PeGr$ (no hermitian structure), $PeQ$ (no hermitian structure,
length~2), $CGr_{0, k}^{0, n}$ and $SCGr_{0, k}^{0, n}$ (no
hermitian structure, lengths $n-k$ and, resp.~$n-k-1$).  The dots
stand for the notions without name or those that require too much
space.  The sign ${\mathfrak a} \; +\hspace{-3.6mm}\supset
{\mathfrak b}$ de\-no\-tes the semidirect sum of Lie
superalgebras, the ideal on the right; Lie superalgebras
${\mathfrak{osp}}_{\alpha}(4|2)$ and ${\mathfrak a}{\mathfrak
b}(3)$ are described via their Cartan matrices, see~\cite{GL3}.

In Table 3:  $m = 2r + 2$ or $m = 2r + 3$, $\varepsilon_1, \ldots,
\varepsilon_r$ and $\delta_1\ldots,\delta_n$ are the standard
bases of the dual spaces to the spaces of diagonal matrices in
${\mathfrak{o}}(m-2)$ and ${\mathfrak{sp}}(n)$, respectively.

In Table 4: $\varepsilon_1, \ldots, \varepsilon_r,
\delta_1\ldots,\delta_n$ is the standard basis of the dual space
to the space of diagonal matrices in ${\mathfrak{gl}}(r|n)$. In
Table 5: $\varepsilon_1, \ldots, \varepsilon_m,
\delta_1\ldots,\delta_{n-m}$ is the standard basis of the dual
space to the space of diagonal matrices in
${\mathfrak{gl}}(m|n-m)$.

\vspace{-3mm}

\begin{table}[th]
\small

\caption{Hermitian symmetric spaces}
\vspace{2mm}

\centering\begin{tabular}{|@{\,}c@{\,}|@{\,}c@{\,}|@{\,}c@{\,}|@{\,}c@{\,}|@{\,}c@{\,}|@{\,}c@{\,}|@{\,}c@{\,}|}
\hline
&&&&&&\\
[-3mm] &\parbox[c]{18mm}{\centering Name\\ of CHSS $X$}
&$X=S_c/G_c$&${\mathfrak{s}}_0=({\mathfrak g}_c)^{{\mathbb C}}$&
${\mathfrak{s}}_{-1}\sim T_0X$&$ (S_c)^*$& \parbox[c]{18mm}{\centering
Name\\ of NCHSS}\\[2mm] \hline
0&$\vphantom{\Big|}{\mathbb C} P^n $ &$SU(n+1)/U(n) $ &$
{\mathfrak{gl}} (n) $& $\id$ &$ SU(1, n)$ &${}^*{\mathbb C}{\mathbb P}
^n$\\
\hline
1&$\vphantom{\Big|}Gr_p^{p+q} $ &$SU(p+q)/S(U(p)\times U(q) $ &$
{\mathfrak{gl}} (p)\oplus {\mathfrak{gl}}(q)$&
$\id\bigodot\id^* $ &$SU(p, q)$ &$ {}^*Gr_p^{p+q}$\\
\hline
2&$\vphantom{\Big|}OGr_n $ &$SO(2n)/U(n) $ &${\mathfrak{gl}}(n) $ &$
\Lambda ^2\id $ &$SO(n, n)$ &$ {}^*OGr_n$\\
\hline
3&$\vphantom{\Big|}LGr_n $ &$Sp(2n)/U(n) $ &${\mathfrak{gl}}(n) $
&$S^2\id$&$Sp(2n; {\mathbb R}) $ &${}^*LGr_n$\\
\hline
4&$\vphantom{\Big|}Q_n $ &$SO(n+2)/(SO(2)\times
SO(n))$&${\mathfrak{co}}(n)$&$\id$&$SO(n, 2)$ &$ {}^*Q_n$\\
\hline
$E_7$&$\vphantom{\Big|}{\mathbb O} P^2$ &$E_7 /(SO(10)\times
U(1))$&${\mathfrak{co}}(10)$&$R(\pi_{5}) $&$E_6^* $&${}^*({\mathbb O} P^2)$\\
\hline
$E_8$&$\vphantom{\Big|}{\mathbb E}$&$E_8 /(E_7 \times U(1)) $
&${\mathfrak{ce}}_7 $ &$R(\pi_{1})$&$E_7 ^* $&${}^*{\mathbb E}$\\
\hline
\end{tabular}
\end{table}

{\bf Occasional isomorphisms:} $Gr_p^{p+q} \cong Gr_q^{p+q}$, $Q_1
\cong {\mathbb C} P ^1$, $Q_2 \cong S^2\times S^2$, $Q_3 \cong LGr_2$,
$Q_4 \cong Gr_2^4$, $OGr_2 \cong LGr_1 \cong {\mathbb C} P ^1$, $OGr_3
\cong Gr_3^4$.

\newpage

\begin{landscape}

\begin{table}[th]
\small
\begin{center}
\caption{Classical superspaces of depth 1}
\vspace{2mm}

\begin{tabular}{|@{\,\;}c@{\,\;}|@{\,\;}c@{\,\;}|@{\,\;}c@{\,\;}|@{\,\;}c@{\,\;}|@{\,\;}c@{\,\;}|@{\,\;}c@{\,\;}|}
\hline
&&&&&\\[-3mm]
${\mathfrak g} $&$ {\mathfrak g}_0$&$ {\mathfrak g}_{-1}$&
Interpretation&\parbox[c]{20mm}{\centering Underlying\\ domain}&
\parbox[c]{24mm}{\centering Name of the\\ superdomain}\\[2mm]
\hline
&&&&&\\[-3mm]
${\mathfrak{sl}}(m| n)$&
$ {\mathfrak{s}}({\mathfrak{gl}}(p| q) \oplus {\mathfrak{gl}}(m-p| n-q))$& $\id\otimes \id^*$&
Supergrassmannian  & $Gr_p^m\times Gr_q^n$&$Gr_{p,  q}^{m,  n}$\\
&&&of the $(p|q)$-dimensional &&\\
&&& subsuperspaces
in ${\mathbb C}^{m| n}$&&\\
${\mathfrak{psl}}(m|m)$&${\mathfrak{ps}}({\mathfrak{gl}}(p| p)\oplus
{\mathfrak{gl}}(m-p|m-p))$&$\id\otimes \id^*$& Same
for $m=n$, $p=q$&$G_p^m\times Gr_p^m$&$Gr_{p,  p}^{m,  m}$\\[1mm]
\hline
&&&&&\\[-3mm]
${\mathfrak{osp}}(m| 2n)$&$ {\mathfrak{cosp}}(m-2| 2n)$& $\id$&
Superquadric of
$(1|0)$-dimen-& $Q_{m-2}$& $Q_{m-2,  n}$\\
 &&& sional isotropic (with &&\\
 &&&
respect to a non-degenerate &&\\
&&& even form) lines in
${\mathbb C}^{m|n}$&&\\[1mm]
 \hline
&&&&&\\[-3mm]
${\mathfrak{osp}}(2m| 2n)$& ${\mathfrak{gl}}(m|n)$ & $\Lambda^2(\id)$& Orthosymplectic Lagrangian
 & ${}^*OGr_m\times LGr_n $&$OLGr_{m,  n}$ \\
&&&supergrassmannian of $(m|m)$-di- &&\\ &&&mensional isotropic (with respect && \\
&&& to a nondegenerate even form)&&\\
 &&&subsuperspaces in ${\mathbb C}^{2m| n}$
&&\\[1mm]
 \hline
&&&&&\\[-3mm]
${\mathfrak{sq}}(n)$&$  {\mathfrak{s}}({\mathfrak{q}}(p)
\oplus {\mathfrak{q}}(n-p))$&$\id\bigodot \id^*$&
Queergrassmannian of  &$Gr_p^n$&$QGr_p^n$\\
${\mathfrak{psq}}(n)$&${\mathfrak{ps}}({\mathfrak{q}}(p) \oplus {\mathfrak{q}}(n-p))$&  &
$\Pi$-symmetric $(p|p)$-dimensional &&\\
&&&subsuperspaces in ${\mathbb C}^{n|n}$&&\\[1mm]
\hline
&&&&&\\[-3mm]
${\mathfrak{pe}} (n)$&${\mathfrak{cpe}}(n-1)$&$ \id$&
Periplectic superquadric of &${\mathbb C} P^{n-1}$&
$PeQ_{n-1}$ \\
${\mathfrak{spe}}(n)$&${\mathfrak{cspe}}(n-1)$& &$(1| 0)$-dimensional
isotropic &&\\
&&&(with respect to a nondegenerate &&\\
&&& odd form) lines in ${\mathbb C} ^{n|n}$ & &\\[1mm]
 \hline
&&&&&\\[-3mm]
${\mathfrak{pe}}(n)$&${\mathfrak{gl}}(p|n-p)$&$\Pi(S^2 (\id))$& Periplectic  Lagrangian
 &$Gr_p^n$& $PeGr_p^n$\\
$({\mathfrak{spe}} (n))$&$
({\mathfrak{sl}}(p|n-p))$& or $\Pi (\Lambda^2 (\id))$& supergrassmannian of $(p|n-p)$-
&& \\
&&& dimensional (and
with a fixed &&\\
&&& volume for ${\mathfrak{spe}}$) subsuperspaces&&\\
&&& in
${\mathbb C}^{n|n}$ isotropic with respect &&\\
&&&  to an odd symmetric &&\\
&&&or skew-symmetric form& &\\ \hline
\end{tabular}
\end{center}\vspace{-1mm}
\end{table}

\end{landscape}

\newpage

\begin{landscape}

\begin{table}[th]

\vspace*{24mm}
\small
\begin{center}
{\bf Table 2 (continued).} Classical superspaces of depth 1

\vspace{2mm}

\begin{tabular}{|@{\,\;}c@{\,\;}|@{\,\;}c@{\,\;}|@{\,\;}c@{\,\;}|@{\,\;}c@{\,\;}|@{\,\;}c@{\,\;}|@{\,\;}c@{\,\;}|}
\hline
&&&&&\\[-3mm]
${\mathfrak g} $&$ {\mathfrak g}_0$&$ {\mathfrak g}_{-1}$&
Interpretation&\parbox[c]{20mm}{\centering Underlying\\ domain}&
\parbox[c]{24mm}{\centering Name of the\\ superdomain}\\[2mm]
\hline
&&&&&\\[-3mm]
${\mathfrak{osp}}_{\alpha}(4|2)$&${\mathfrak{cosp}}
 (2|2)\simeq {\mathfrak{gl}}(2|1)$ & $\id$&$\cdots$& ${\mathbb C}
P^1\times {\mathbb C} P^1$ &$\cdots$\\
${}={\mathfrak{osp}} (4| 2)_{\alpha}$&&&&&\\[1mm]
\hline
&&&&&\\[-3mm] ${\mathfrak{ab}}(3)$&${\mathfrak{cosp}}(2|4)$&$
L_{3\varepsilon_{1}}$ &$\cdots$& ${\mathbb C} P^1\times Q_5$
&$\cdots$\\[1mm]
\hline
${\mathfrak{vect}} (0| n)$&${\mathfrak{vect}}(0| n-k)\; +\hspace{-3.6mm}\supset
 {\mathfrak{gl}}(k; \Lambda (n-k))$
&$\Lambda (k)\otimes \Pi (\id)$&Curved supergrassmanian&$\cdots$ &$CGr_{0,
k}^{0,  n}$\\
 &&& of $(0|1)$-dimensional && \\
&&& subsupermanifolds in
${\mathbb C}^{0|n}$&&\\
${\mathfrak{svect}} (0| n)$&${\mathfrak{vect}}(0| n-k)\;
+\hspace{-3.6mm}\supset {\mathfrak{sl}}(k; \Lambda (n-k))$& $\Pi
(Vol)$ if $k=1$ &Same with volume elements &$\cdots$&$SCGr_{0, k}^{0,
n}$\\
&&&preserved in the sub- and ambient &&\\
&&& supermanifolds &&\\[1mm]
\hline
&&&&&\\[-3mm]
${\mathfrak{h}}(0| m)$&${\mathfrak{h}} (0| m-2)\; +\hspace{-3.6mm}
\supset \Lambda (m-2)\cdot z$&& Curved
superquadric&$\cdots$&$CQ_{m-2,  0}$ \\
&&& of $(0|1)$-dimensional &&\\
&&& subsupermanifolds in ${\mathbb C}^{0|m}$&&\\
$ {\mathfrak{h}}^o(m)$&$ {\mathfrak{h}}^o(m-2)\; +\hspace{-3.6mm}\supset
\Lambda (m-2)\cdot z$ &$ \Pi (\id)$& isotropic (with respect to a
(partly)&&\\
&&& split symmetric form)&&\\ \hline
\end{tabular}
\end{center}
\end{table}

\end{landscape}

\newpage

\begin{table}
\small
\begin{minipage}[t]{0.42\linewidth}
\caption{}
\vspace{2mm}
\begin{tabular}{|@{\;}c@{\;}|@{\;}c@{\;}|@{\;}c@{\;}|@{\;}c@{\;}|}
\hline
$\vphantom{\Big|}r$&$n$&$H^{2,  2}_{({\mathfrak{g}} _{-1},
 \widehat{{\mathfrak{g}}} _{0})_{*}}$&
$S^  2({\mathfrak{g}} _{-1})$\\
\hline
$\vphantom{\Big|}0$&$1$&---&$0, \ \mbox{if}\ m=2 $\\
 & & &$\vphantom{\Big|}0,  \delta _1\ \mbox{if}\ m=3 $\\
\hline
$\vphantom{\Big|}0$&$ \geq 2$ &$ 2\delta _1+2\delta _2$&$0,  \delta _1+ \delta _2 $\\
\hline
$\vphantom{\Big|}1$&$ 1$&$\varepsilon _1+\delta _ 1$&$ 0,   2\varepsilon _1 $\\
\hline
$\vphantom{\Big|}1$&$  \geq 2$&$2\varepsilon _1+\delta _ 1+\delta _ 2$&$0,   2\varepsilon _1 $\\
\hline
$\vphantom{\Big|} \geq 2$&$ \geq 1$&$2\varepsilon _1+2\varepsilon _2$&$0,   2\varepsilon _1$\\
\hline
\end{tabular}
\end{minipage}
\hfill
\begin{minipage}[t]{0.52\linewidth}
\caption{}
\vspace{2mm}
\begin{tabular}{|@{\;}c@{\;}|@{\;}c@{\;}|@{\;}c@{\;}|@{\;}c@{\;}|}
\hline
$\vphantom{\Big|}r$&$n$&$H^{1,  2}_{({\mathfrak{g}} _{-1},
 {{\mathfrak{g}}}_{0})_{*}}$&
$H^{2,  2}_{({\mathfrak{g}}_{-1},  {{\mathfrak{g}}} _{0})_{*}}$\\
\hline
$\vphantom{\Big|} \geq 5$&$0$&$\varepsilon _1+\varepsilon _2 -\varepsilon _{r-2}-\varepsilon _{r-1}-2\varepsilon _{r}$
&---\\
\hline
$\vphantom{\Big|}0$&$ \geq 3$ &$ 2\delta _1-\delta _{n-1}-3\delta _{n }$&---\\
\hline
$\vphantom{\Big|}1$&$ \geq 3$&$\varepsilon _1+\delta _1-\delta _{n-1}-3\delta _{n }$&---\\
\hline
$\vphantom{\Big|} \geq 5$&$ 1$&$\varepsilon _1+\varepsilon _2 -\varepsilon _{r}-3\delta _{n }$
&---\\
\hline
$\vphantom{\Big|} \geq 2$&$ \geq 2$&$\varepsilon _1+\varepsilon _2 -\delta _{n-1}-3\delta _{n }$&
---\\
\hline
$\vphantom{\Big|} 2$&$1$&---&$-\varepsilon _2-3\delta _1$\\
\hline
\end{tabular}
\end{minipage}
\end{table}

\vspace{-7mm}

\begin{table}[th]\small
\begin{center}
\caption{}
\vspace{2mm}
\begin{tabular}{|@{\;}c@{\;}|@{\;}c@{\;}|@{\;}c@{\;}|}
\hline
$\vphantom{\Big|}m$&$n-m$&$H^{1,  2}_{({\mathfrak{g}} _{-1},
  {\mathfrak{g}} _{0})_{*}}$\\
\hline
$\vphantom{\Big|}0$&$ \geq 4$&$\delta_1+\delta_2-2\delta_{n-1}-2\delta_n;
 \quad
-\delta_{n-1}-\delta_n$\\
\hline
$\vphantom{\Big|} \geq 1$&$ \geq 2$&$2\varepsilon_1-2\delta_{n-m-1}-2\delta_{n-m};
  \quad
-\delta_{n-m-1}-\delta_{n-m}$\\
\hline
$\vphantom{\Big|}2$&$1$&$2\varepsilon_1-3\varepsilon_2-\delta_{1};  \quad
-\varepsilon_2-\delta_1$\\
\hline
$\vphantom{\Big|} \geq 3$&$0$&$2\varepsilon_1-4\varepsilon_m;  \quad 2\varepsilon_1-2\varepsilon_{m-1}-2\varepsilon_m;  \quad
-2\varepsilon_m$\\
\hline
$\vphantom{\Big|} \geq 3$&$1$&$2\varepsilon_1-\varepsilon_{m-1}-\varepsilon_m-2\delta_1;  \quad
2\varepsilon_1-3\varepsilon_{m}-\delta_1$\\
\hline\end{tabular}
\end{center}
\vspace{-3mm}

\end{table}
\section{Structure functions for exceptional superdomains}

In this section we set $Y_i=X_i^+$, $X_i=X_i^-$.

{\bf \mathversion{bold}${\mathfrak{g}}={\mathfrak{osp}}_{\alpha}(4|2)$.} There are five types of Cartan matrices obtained from each other via odd reflections, see \cite{GL3}, but we only consider realization with the Cartan matrices
\[
1)\; \begin{pmatrix} 0&1&-1-\alpha\\ -1&0&-\alpha\\
-1-\alpha&\alpha&0\end{pmatrix}
 \quad \text{
and}
 \quad
 2)\; \begin{pmatrix} 2&-1&0\\
-1&0&-\alpha\\
0&-1&2\end{pmatrix}\] 
because the classical superdomains corresponding to other matrices are the same as the ones obtained from these Cartan matrices. 

\begin{theorem} The structure functions are only of order $2$.

1) For the parabolic subalgebra generated by $X^+_{1}$, $X^{\pm}_2$,
$X^{\pm}_3$, set $X^{-}_{4}=[X^{-}_1, X^{-}_2]$, $X^{-}_{5}=[X^{-}_1,
X^{-}_3]$, $X^{-}_{6}=[X^{-}_2, X^{-}_3]$, $X^{-}_{7}=[X^{-}_1,
[X^{-}_2, X^{-}_3]]$.

The ${\mathfrak{g}}_{0}={\mathfrak{gl}}(1|2)$-module of structure functions is irreducible,
the highest weight vector (a representative of the cohomology class)
is odd and its highest weight (with respect to Borel subalgebra given by
$X^{+}_2$, $X^{+}_3$) in basis $H_{1}$, $H_{2}$, $H_3$ are as follows:
\begin{center}
\begin{tabular}{|c|}
\hline
$-\alpha(\alpha+1)H_1dY_4dY_7+\alpha^2H_2dY_4dY_7+(1+\alpha)X_2dY_4dY_5+
X_6dY_1dY_4\vphantom{\Big|}$\\
\hline
$\left(-1, -\frac{4}{\alpha}-1, 1\right)\vphantom{\Big|}$\\
\hline
\end{tabular}
\end{center}

2) For the parabolic subalgebra generated by $X^+_{1}$, $X^{\pm}_2$,
$X^{\pm}_3$, set $X^{-}_{4}=[X^{-}_1, X^{-}_2]$, $X^{-}_{5}=[X^{-}_2,
X^{-}_3]$, $X^{-}_{6}=[X^{-}_3, [X^{-}_1, X^{-}_2]]$,
$X^{-}_{7}=[[X^{-}_1, X^{-}_2], [X^{-}_2, X^{-}_3]]$.

The ${\mathfrak{g}}_{0}={\mathfrak{gl}}(1|2)$-module of structure functions is irreducible,
the highest weight vector (a representative of the cohomology class)
is odd and its heighest weight (with respect to Borel subalgebra given
by $X^{+}_1$, $X^{+}_2$) in basis $H_{1}$, $H_2$, $H_3$ are as follows:
\begin{center}
\begin{tabular}{|c|}
\hline
$(2+\alpha)H_1dY_1dY_6+2H_2dY_1dY_6+\alpha H_3dY_1dY_6+
2(1+\alpha)Y_2dY_1dY_7\vphantom{\Big|}$\\
\hline
$(-1, -1-\alpha, 1)\vphantom{\Big|}$\\
\hline
\end{tabular}
\end{center}

\end{theorem}

{\bf \mathversion{bold}${\mathfrak{g}}={\mathfrak{ab}}_3$.} We
consider realization with the Cartan matrix
\[\begin{pmatrix}
2&-1&0&0\\ -3&0&1&0\\
0&-1&2&-2\\ 0&0&-1&2\end{pmatrix} \] 


\begin{theorem} The structure functions are only of order $1$.  The
${\mathfrak{g}}_{0}={\mathfrak{osp}}(2|4)$-module of structure functions is reducible, so we describe it in terms of the highest weight vectors (representatives of the cohomology classes) with respect to $({\mathfrak{g}}_{0})_{\bar 0}$) in the following table (where $p$ is parity, ${\rm deg}$ is the degree ($\deg x_1=\dots=\deg x_4=1$; in general this \lq\lq$\deg$" is more convenient than weight which may depend on a complex parameter, e.g., as for ${\mathfrak{osp}}_{\alpha}(2|4)$); ${\mathfrak{sp}}$ is the weight wrt ${\mathfrak{sp}}(4)$ (it does not matter how it is embedded; the weight is given for those who wonder how the dimensions were computed) and the  operator that grades ${\mathfrak{g}}$):

\begin{center}
\begin{tabular}{|l|l|l|l|l|l|l|}
\hline
\#&weight&${\mathfrak{sp}}$&{\rm deg}&vector&$p$&$\dim$\\
\hline
$1$&$3, -2, 1, 1$&$2, 1$&$1, -1, 1, 1$&$Y_{11}dY_{1}dY_{15}$&$1$&$16$\\
\hline
$2$&$2, -2, 0, 1$&$1, 1$&$1, 0, 1,
1$&$Y_{5}dY_{1}dY_{12}+Y_{8}dY_{1}dY_{15}$&$0$&$5$\\
\hline
$3$&$2, -1, 2, 0$&$2, 0$&$1, 0, 2, 1$&$Y_{5}dY_{1}dY_{15}$&$0$&$10$\\
\hline
$4$&$2, -1, 0, 2$&$2, 2$&$1, 0, 2, 2$&$Y_{11}dY_{1}dY_{18}$&$0$&$14$\\
\hline
$5$&$2, 0, 2, 1$&$3, 1$&$1, 0, 3, 2$&$Y_{11}dY_{15}dY_{15}$&$0$&$35$\\
\hline
$6$&$1, -1, 1, 0$&$1, 0$&$1, 1, 2, 1$&$Y_{1}dY_{1}dY_{15}$&$1$&$4$\\
\hline
$7$&$1, 0, 1, 1$&$2, 1$&$1, 1, 3, 2$&a: $Y_{5}dY_{1}dY_{18}$&$1$&$16$\\
&&&&b: $Y_{5}dY_{12}dY_{15}+Y_{8}dY_{15}dY_{15}$&&$16$\\
\hline
$8$&$1, 1, 3, 0$&$3, 0$&$1, 1, 4, 2$&$Y_{5}dY_{15}dY_{15}$&$1$&$20$\\
\hline
$9$&$1, 1, 1, 2$&$3, 2$&$1, 1, 4, 3$&$Y_{1}dY_{11}dY_{18}$&$1$&$40$\\
$10$&$0, 0, 0, 1$&$1, 1$&$1, 2, 3, 2$&$Y_{1}dY_{1}dY_{18}$&$0$&$5$\\
\hline
$11$&$0, 1, 2, 0$&$2, 0$&$1, 2, 4, 2$&a: $Y_{1}dY_{1}dY_{15}$&$0$&$10$\\
&&&&b: $-2Y_{5}dY_{8}dY_{18}-4Y_{5}dY_{12}dY_{17}$&&$10$\\
&&&&$+Y_{5}dY_{15}dY_{16}$&&\\
\hline
$12$&$0, 1, 0, 2$&$2, 2$&$1, 2, 4, 3$&$Y_{5}dY_{12}dY_{18}+
Y_{8}dY_{15}dY_{18}$&$0$&$14$\\
\hline
$13$&$0, 2, 2, 1$&$3, 1$&$1, 2, 5, 3$&$Y_{5}dY_{15}dY_{18}$&$0$&$35$\\
\hline
$14$&$-1, 1, 1, 0$&$1, 0$&$1, 3, 4, 2$&$-2Y_{1}dY_{8}dY_{18}-4
Y_{1}dY_{12}dY_{17}$&$1$&$4$\\
&&&&$+Y_{1}dY_{15}dY_{16}$&&\\
\hline
$15$&$-1, 2, 1, 1$&$2, 1$&$1, 3, 5, 3$&$Y_{1}dY_{15}dY_{18}$&$1$&$16$\\
\hline
$16$&$-1, 3, 3, 0$&$3, 0$&$1, 3, 6, 3$&$Y_{5}dY_{17}dY_{18}$&$1$&$20$\\
\hline
$17$&$-2, 3, 2, 0$&$2, 0$&$1, 4, 6, 3$&$Y_{1}dY_{17}dY_{18}$&$0$&$10$\\
\hline\end{tabular}
\end{center}

Irreducible ${\mathfrak{osp}}(2|4)$-submodules of $H^{1, 2}$:
\begin{gather*}
    A=[1]\oplus [2]\oplus[4]\oplus[5]\oplus [7:
2a-3b]\oplus [9]\oplus[12]  \quad  (\dim =68|72),\\
B=[11:a-b]\oplus[14]\oplus[16]\oplus [17]  \quad  (\dim =20|24).
\end{gather*}
The quotient of $H^{1, 2}$ modulo $A\oplus B$ is an irreducible
${\mathfrak{osp}}(2|4)$-module.

The maximal parabolic
subalgebra corresponds to the first Chevalley generator;
\begin{gather*}
    X^{-}_{5}=[X^{-}_1, X^{-}_2],  \quad X^{-}_{6}=[X^{-}_2, X^{-}_3],
 \quad X^{-}_{7}=[X^{-}_3, X^{-}_4], \\
X^{-}_{8}=[X^{-}_3, [X^{-}_1,
X^{-}_2]],  \quad X^{-}_{9}=[X^{-}_3, [X^{-}_3,
X^{-}_4]],  \quad X^{-}_{10}=[X^{-}_4, [X^{-}_2,
X^{-}_3]], \\
X^{-}_{11}=[[X^{-}_1, X^{-}_2], [X^{-}_2, X^{-}_3]],  \quad
X^{-}_{12}=[[X^{-}_1, X^{-}_2], [X^{-}_3, X^{-}_4]],\\
X^{-}_{13}=[[X^{-}_2, X^{-}_3], [X^{-}_3, X^{-}_4]],
 \quad X^{-}_{14}=[[X^{-}_1, X^{-}_2], [X^{-}_4, [X^{-}_2, X^{-}_3]]],\\
X^{-}_{15}=[[X^{-}_3, X^{-}_4], [X^{-}_3, [X^{-}_1, X^{-}_2]]],
\\
 X^{-}_{16}=[[X^{-}_3, [X^{-}_1, X^{-}_2]], [X^{-}_4, [X^{-}_2,
X^{-}_3]]],\\
X^{-}_{17}=[[X^{-}_3, [X^{-}_3, X^{-}_4]], [[X^{-}_1, X^{-}_2], [X^{-}_2,
X^{-}_3]]], \\
X^{-}_{18}=[[[X^{-}_1, X^{-}_2], [X^{-}_3, X^{-}_4]], [[X^{-}_2,
X^{-}_3], [X^{-}_3, X^{-}_4]]].
\end{gather*}
\end{theorem}

\section{Structure functions for the \lq\lq odd Penrous" tensor}

All ${\mathbb Z}$-gradings of depth 1 of ${\mathfrak g} = {\mathfrak{psq}}(n)$ are of the form
${\mathfrak g}_{-1} \oplus {\mathfrak g}_0 \oplus {\mathfrak g}_{1}$, where ${\mathfrak g}_0 = {\mathfrak {ps}} ({\mathfrak q}(p) \oplus
{\mathfrak q}(n-p))$ for $p > 0$, and 
${\mathfrak g}_1 \cong {\mathfrak g}_{-1}^*$,  as ${\mathfrak g}_0$-modules, where ${\mathfrak g}_{-1}$ is either one of the two irreducible
${\mathfrak g}_0$-modules in $V(p|p)^* \otimes V(n-p|n-p)$. Explicitly:
$$
{\mathfrak g}_{-1} = \Span ((x \pm \Pi(x)) \otimes (y \pm \Pi(y)) ), \text{ where }
x\in V(p|p)^*, \quad y\in V(n-p|n-p).
$$
Let $ \varepsilon_1, \ldots ,
 \varepsilon_p$ and $\delta_1, \ldots , \delta_{n-p}$ be the standard bases of the dual  spaces to the spaces of diagonal matrices in ${\mathfrak q}(p)$ and ${\mathfrak q}(n-p)$, 
respectively.

\begin{theorem} 1) $({\mathfrak g}_{-1}, {\mathfrak g}_0)_* = {\mathfrak g}$. 
 
2) all structure functions are of order 1; they split into the direct sum of two irreducible  ${\mathfrak g}_0$-submodules with highest weights 
$2 \varepsilon_1 -  \varepsilon_p + \delta_1 - 2\delta_{n-p}$ and
$\varepsilon_1 - \delta_{n-p}$. \end{theorem}

\appendix

\section*{Appendix}

\section{Background}

{\bf Linear algebra in superspaces.  Generalities.} A {\it superspace}
is a ${\mathbb Z} /2$-graded space; for a superspace $V=V_{\bar
0}\oplus V_{\bar 1}$ denote by $\Pi (V)$ another copy of the same
superspace: with the shifted parity, i.e., $(\Pi(V))_{\bar i}= V_{\bar
i+\bar 1}$.  The {\it superdimension} of $V$ is $\dim
V=p+q\varepsilon$, where $\varepsilon^2=1$ and $p=\dim V_{\bar 0}$,
$q=\dim V_{\bar 1}$.  (Usually $\dim V$ is shorthanded as a pair $(p,
q)$ or $p|q$; with the help of $\varepsilon$ the fact that $\dim
V\otimes W=\dim V\cdot \dim W$ becomes lucid.)

A {\it superalgebra} is a superspace $A$ with an
even multiplication map $m: A\otimes A\longrightarrow A$.

A superspace structure in $V$ induces the superspace structure in the
space $\End (V)$.  A~{\it basis of a superspace} always consists of
{\it homogeneous} vectors; let $\Par=(p_1, \dots, p_{\dim V})$ be an
ordered collection of their parities.  We call $\Par$ the {\it format}
of the basis of $V$.  A~square {\it supermatrix} of format (size)
$\Par$ is a $\dim V\times \dim V$ matrix whose $i$th row and $i$th
column are of the same parity $p_i$.  The matrix unit $E_{ij}$ is
supposed to be of parity $p_i+p_j$ and the bracket of supermatrices
(of the same format) is defined via Sign Rule: {\it if something of
parity $p$ moves past something of parity $q$ the sign $(-1)^{pq}$
accrues; the formulas defined on homogeneous elements are extended to
arbitrary ones via linearity}.

Examples of application of Sign Rule: setting $[X,
Y]=XY-(-1)^{p(X)p(Y)}YX$ we get the notion of the supercommutator and
the ensuing notions of the supercommutative superalgebra and the Lie
superalgebra (that in addition to superskew-commutativity satisfies
the super Jacobi identity, i.e., the Jacobi identity amended with the
Sign Rule).  The {\it superderivation} of a superalgebra $A$ is a
linear map $D: A\longrightarrow A$ such that satisfies the Leibniz
rule (and Sign rule)
\[
D(ab)=D(a)b+(-1)^{p(D)p(a)}aD(b).
\]

Usually, $\Par$ is of the form $(\bar 0 , \dots, \bar 0 , \bar 1 ,
\dots, \bar 1)$.  Such a format is called {\it standard}.  In this
paper we can do without nonstandard formats.  But they are vital in
various questions related with the study of distinct systems of simple
roots that the reader might be interested in.

The {\it general linear} Lie superalgebra of all supermatrices of size
$\Par$ is denoted by ${\mathfrak{gl}}(\Par)$; usually,
${\mathfrak{gl}}(\bar 0, \dots, \bar 0, \bar 1, \dots, \bar 1)$ is
abbreviated to ${\mathfrak{gl}}(\dim V_{\bar 0}|\dim V_{\bar 1})$.
Any matrix from ${\mathfrak{gl}}(\Par)$ can be expressed as the sum of
its even and odd parts; in the standard format this is the block
expression:
\[
\begin{pmatrix}A&B\\ C&D\end{pmatrix}=\begin{pmatrix}A&0\\
0&D\end{pmatrix}+\begin{pmatrix}0&B\\ C&0\end{pmatrix},\quad
p\left(\begin{pmatrix}A&0\\
0&D\end{pmatrix}\right)=\bar 0, \quad p\left(\begin{pmatrix}0&B\\
C&0\end{pmatrix}\right)=\bar 1.
\]

The {\it supertrace} is the map ${\mathfrak{gl}} (\Par)\longrightarrow
{\mathbb C}$, $(A_{ij})\mapsto \sum (-1)^{p_{i}}A_{ii}$.  Since $\str
[x, y]=0$, the space of supertraceless matrices constitutes the {\it
special linear} Lie subsuperalgebra ${\mathfrak{sl}}(\Par)$.

There are, however, two super versions of ${\mathfrak{gl}}(n)$, not
one.  The other version is called the {\it queer} Lie superalgebra and
is defined as the one that preserves the complex structure given by an
{\it odd} operator $J$, i.e., is the centralizer $C(J)$ of $J$:
\[
{\mathfrak{q}}(n)=C(J)=\{X\in{\mathfrak{gl}}(n|n): [X, J]=0 \},  \quad
\mbox{where}  \quad J^2=-\id.
\]
It is clear that by a change of basis we can reduce $J$ to the form
$J_{2n}=\begin{pmatrix}0&1_n\\ -1&0\end{pmatrix}$.  In the standard
format we have
\[
{\mathfrak{q}}(n)=\left \{\begin{pmatrix}A&B\\ B&A\end{pmatrix}\right\}.
\]
On ${\mathfrak{q}}(n)$, the {\it queer trace} is defined: ${\mathfrak q}tr:
\begin{pmatrix}A&B\\
B&A\end{pmatrix}\mapsto
\tr B$. Denote by ${\mathfrak{sq}}(n)$ the Lie superalgebra of {\it queertraceless}
matrices.

Observe that the identity representations of ${\mathfrak{q}}$ and
${\mathfrak{sq}}$ in $V$, though irreducible in super sence, are not
irreducible in the nongraded sence: take homogeneous linearly
independent vectors $v_1, \dots , v_n$ from $V$; then $\Span
(v_1+J(v_1), \dots , v_n+J(v_n))$ is an invariant subspace of $V$,
which is not a subsuperspace, singled out by a $\Pi$-symmetry.

We will stick to the following terminology, cf.~\cite{L1,L2}.
The representation of a superalgebra $A$ in the superspace $V$ is
irreducible of {\it general type} or just of {\it $G$-type} if it does
not contain homogeneous (with respect to parity) subrepresentations
distinct from $0$ and $V$ itself, otherwise it is called {\it
irreducible of $Q$-type}.  Thus, an irreducible representation of
$Q$-type has no invariant sub{\it super}space but {\it has} a
nontrivial invariant subspace.

So, there are two types of irreducible representations: those that do
not contain any nontrivial subrepresentations (called of {\it general
type} or of type $G$) and those that contain {\it in}homogeneous
invariant subspaces (called them of {\it type} $Q$).  If $V$ is of
finite dimension, then in the first case its centralizer, as of
$A$-module, is isomorphic to ${\mathfrak{gl}}(1)$, whereas in the second case to
${\mathfrak{q}}(1)$.

Let $V_{1}$ and $V_{2}$ be finite dimensional irreducible modules over
$A_{1}$ and $A_{2}$, respectively.  Then $V_{1}\otimes V_{2}$ is an
irreducible $A_{1}\otimes A_{2}$-module except for the case when both
$V_{1}$ and $V_{2}$ are of type $Q$.  In the latter case, the
centralizer of the $A_{1}\otimes A_{2}$-module $V_{1}\otimes V_{2}$ is
isomorphic to $Cl_{2}$, the Clifford superalgebra with 2 generators.

If $e\in Cl_{2}$ is a minimal idempotent, then $e(V_{1}\otimes V_{2})$
is an irreducible $A_{1}\otimes A_{2}$-module of type $G$ that we will
denote by $V_{1}\bigodot V_{2}$, see Tables~1 and~2.

More generally, we can consider matrices with the elements from a
supercommutative superalgebra $\Lambda$.  Then the parity of the
matrix with only one nonzero $i, j$-th element $X_{i, j}\in \Lambda$
is equal to $p_{i}+p_j+p(X_{i, j})$.

{\bf The berezinian and the module of volume forms.}  On $GL(p|q;
\Lambda)$, the group of even $p|q\times p|q$ invertible matrices with
elements from a supercommutative superalgebra $\Lambda$, a
multiplicative function~--- an analog of determinant~--- is defined.
In honor of F~Berezin Leites baptized it {\it berezinian}, cf.~\cite{L0}.
 Its explicit expression in the standard format is
\[
\ber\begin{pmatrix}A&B\\
C&D\end{pmatrix}=\det\left(A-BD^{-1}C\right)\det D^{-1}
\]
or
\[
\ber^{-1}\begin{pmatrix}A&B\\
C&D\end{pmatrix}=\det\left(D-CA^{-1}B\right)\det A^{-1}.
\]
The berezinian is a rational function and this is a reason why the
structure of the algebra of invariant polynomials on
${\mathfrak{gl}}(p|q)$ is much more complicated than that for the
Lie algeb\-ra~${\mathfrak{gl}}(n)$.

Clearly, the derivative of the berezinian is supertrace and the
relation between them is as expected: $\ber X=\exp\str \log X$.

The one-dimensional representation $\Vol(V)$ of $GL(V; \Lambda)$
corresponding to $\ber$ and at the same time to the representation
$\str$ of ${\mathfrak{gl}}(V)$ is called the space of {\it volume forms}.  It
can be only realized in the space of tensors as a quotient module:
recall that for ${\mathfrak{gl}}(V)$ there is no complete reducibility, cf.~\cite{L3}.

{\bf The odd analog of berezinian.}  On the group $GQ(n; \Lambda)$
of invertible even matrices from $Q(n; \Lambda)$, the berezinian is
identically equal to 1.  Instead, on $GQ(n; \Lambda)$ there is
defined its own {\it queer determinant}
\[
{\mathfrak q}et\begin{pmatrix}A&B\\
B&A\end{pmatrix}=\sum\limits_{i \geq 0}\frac{1}{2i+1}\tr\left(A^{-1}B\right)^{2i+1}.
\]
This strange function is $GQ(n; \Lambda)$-invariant and additive,
i.e., ${\mathfrak q}et XY={\mathfrak q}et X+{\mathfrak q}et Y$, cf.~\cite{BL}.

{\bf Superalgebras that preserve bilinear forms: two types}.  To the
linear map $F: V\longrightarrow W$ of superspaces there corresponds the dual map
$F^*: W^*\longrightarrow V^*$ of the dual superspaces; if $A$ is the supermatrix
corresponding to $F$ in a basis of the format $\Par$, then, in the
dual basis, to $F^*$ the {\it supertransposed} matrix $A^{st}$
corresponds:
\[
(A^{st})_{ij}=(-1)^{(p_{i}+p_{j})(p_{i}+p(A))}A_{ji}.
\]

The supermatrices $X\in{\mathfrak{gl}}(\Par)$ such that
\[
X^{st}B+(-1)^{p(X)p(B)}BX=0 \quad\mbox{ for an homogeneous matrix
$B\in{\mathfrak{gl}}(\Par)$}
\]
constitute the Lie superalgebra ${\mathfrak{aut}} (B)$ that preserves the
bilinear form on $V$ with matrix~$B$.  Most popular is the
nondegenerate supersymmetric form whose matrix in the standard format
is the canonical form $B_{\rm ev}$ or $B'_{\rm ev}$:
\[
B_{\rm ev}(m|2n)= \begin{pmatrix} 1_m&0\\
0&J_{2n}
\end{pmatrix}, \mbox{where}  \quad
J_{2n}=\begin{pmatrix}0&1_n\\-1_n&0\end{pmatrix},
\]
or
\[
B'_{\rm ev}(m|2n)= \begin{pmatrix} \antidiag (1, \dots , 1)&0\\
0&J_{2n}
\end{pmatrix}.
\]
The usual notation for ${\mathfrak{aut}} (B_{\rm ev}(m|2n))$ is ${\mathfrak{osp}}(m|2n)$ or
${\mathfrak{osp}}^{sy}(m|2n)$.

Recall that the ``upsetting'' map $u: \Bil(V, W)\longrightarrow \Bil(W, V)$ becomes
for $V=W$ an involution $u: B\mapsto B^u$ which on matrices acts as
follows:
\[
B=\begin{pmatrix}
B_{11}&B_{12}\\
B_{21}&B_{22}
\end{pmatrix}\mapsto B^u=\begin{pmatrix}
B_{11}^t&(-1)^{p(B)}B_{21}^t\\
(-1)^{p(B)}B_{12}^t&B_{22}^t
\end{pmatrix}.
\]
The forms $B$ such that $B=B^u$ are called {\it supersymmetric} and
{\it superskew-symmetric} if $B=-B^u$.  The passage from $V$ to $\Pi (V)$
identifies the space of supersymmetric forms on $V$ with that
superskew-symmetric ones on $\Pi (V)$, preserved by the
``symplectico-orthogonal" Lie superalgebra ${\mathfrak{osp}}^{sk}(m|2n)$ which is
isomorphic to ${\mathfrak{osp}}^{sy}(m|2n)$ but has a different matrix
realization.  We never use notation ${\mathfrak{sp}}'{\mathfrak{o}} (2n|m)$ in order not to
confuse with the special Poisson superalgebra.

In the standard format the matrix realizations of these algebras
are:
\[
{\mathfrak{osp}} (m|2n)=\left\{\left (\begin{matrix} E&Y&X^t\\
X&A&B\\
-Y^t&C&-A^t\end{matrix} \right)\right\},\quad {\mathfrak{osp}}^{sk}(m|2n)=
\left\{\left(\begin{matrix} A&B&X\\
C&-A^t&Y^t\\
Y&-X^t&E\end{matrix} \right)\right\},
\]
where $\left(\begin{matrix} A&B\\
C&-A^t\end{matrix} \right)\in {\mathfrak{sp}}(2n)$, $E\in{\mathfrak{o}}(m)$
and  ${}^t$ is the usual transposition.

A nondegenerate supersymmetric odd bilinear form $B_{\rm odd}(n|n)$ can be
reduced to the canonical form whose matrix in the standard format is
$J_{2n}$.  A canonical form of the superskew odd nondegenerate form in
the standard format is $\Pi_{2n}=\begin{pmatrix}
0&1_n\\1_n&0\end{pmatrix}$.  The usual notation for ${\mathfrak{aut}}
(B_{\rm odd}(\Par))$ is ${\mathfrak{pe}}(\Par)$.  The passage from $V$ to $\Pi (V)$
sends the supersymmetric forms to superskew-symmetric ones and
establishes an isomorphism ${\mathfrak{pe}}^{sy}(\Par)\cong{\mathfrak{pe}}^{sk}(\Par)$.
This Lie superalgebra is called, as A~Weil suggested to Leites, {\it
periplectic}.  In the standard format these superalgebras are
shorthanded as in the following formula, where their matrix
realizations is also given:
\begin{gather*}
{\mathfrak{pe}} ^{sy} (n)=\left\{\begin{pmatrix} A&B\\
C&-A^t\end{pmatrix}, \ \mbox{where}\ B=-B^t,\
C=C^t\right\};\\
{\mathfrak{pe}}^{sk}(n)=\left\{\begin{pmatrix}A&B\\ C&-A^t\end{pmatrix}, \
\mbox{where}\ B=B^t, \ C=-C^t\right\}.
\end{gather*}

The {\it special periplectic} superalgebra is ${\mathfrak{spe}}(n)=\{X\in{\mathfrak{pe}}(n): \str
X=0\}$.

Observe that though the Lie superalgebras ${\mathfrak{osp}}^{sy} (m|2n)$ and
${\mathfrak{pe}} ^{sk} (2n|m)$, as well as ${\mathfrak{pe}} ^{sy} (n)$ and
${\mathfrak{pe}} ^{sk} (n)$,
are isomorphic, the difference between them is sometimes crucial.

{\bf  Projectivization.} If ${\mathfrak{s}}$ is a Lie algebra of scalar
matrices, and ${\mathfrak{g}}\subset {\mathfrak{gl}} (n|n)$ is a Lie subsuperalgebra
containing ${\mathfrak{s}}$, then the {\it projective} Lie superalgebra of type
${\mathfrak{g}}$ is ${\mathfrak{pg}}= {\mathfrak{g}}/{\mathfrak{s}}$.

Projectivization sometimes leads to new Lie superalgebras, for
example: ${\mathfrak{pgl}} (n|n)$, ${\mathfrak{psl}} (n|n)$, ${\mathfrak{pq}} (n)$,
${\mathfrak{psq}} (n)$;
whereas ${\mathfrak{pgl}} (p|q)\cong {\mathfrak{sl}} (p|q)$ if $p\neq q$.

\section{Certain constructions with the point functor}
The point functor is well-known in algebraic geometry since at least
1953~\cite{Wi}.  The advertising of ringed spaces with nilpotents in
the structure sheaf that followed the discovery of supersymmetries
caused many mathematicians and physicists to realize the usefulness of
the language of points; most interesting are numerous ideas due to
Witten (for some of them see~\cite{W}), for their clarifications, and
further developments and references see \cite{M,D}.
F~A~Berezin~\cite{B} was the first who applied the point functor to
study Lie superalgebras.

All superalgebras and modules are supposed to be
finite dimensional over ${\mathbb C} $.

{\bf What a Lie superalgebra is.} Lie superalgebras had appeared
in topology in 1930's or earlier.  So when somebody offers a
``better than usual'' definition of a notion which seemed to have
been established about 70 year ago this might look strange, to say
the least.  Nevertheless, the answer to the question ``what is a
Lie superalgebra?''  is still not a~common knowledge.  Indeed, the
naive definition (``apply the Sign Rule to the definition of the
Lie algebra'') is manifestly inadequate for considering the
(singular) supervarieties of deformations and applying
representation theory to mathematical physics, for example, in the
study of the coadjoint representation of the Lie supergroup which
can act on a~supermanifold but never on a superspace (an object
from another category).  So, to deform Lie superalgebras and apply
group-theoretical methods in ``super'' setting, we must be able to
recover a supermanifold from a superspace, and vice versa.

A proper definition of Lie superalgebras is as follows, cf.~\cite{L1,D}.  The {\it Lie superalgebra} in the category of
supermanifolds corresponding to the ``naive'' Lie superalgebra $L=
L_{\bar 0} \oplus L_{\bar 1}$ is a linear supermanifold ${\mathcal L}=(L_{\bar 0},
{\mathcal O})$, where the sheaf of functions ${\mathcal O}$ consists of functions on
$L_{\bar 0}$ with values in the Grassmann superalgebra on $L_{\bar 1}^*$;
this supermanifold should be such that for ``any'' (say, finitely
generated, or from some other appropriate category) supercommutative
superalgebra $C$, the space ${\mathcal L}(C)=\Hom (\Spec C, {\mathcal L})$, called {\it
the space of $C$-points of} ${\mathcal L}$, is a Lie algebra and the
correspondence $C\longrightarrow {\mathcal L}(C)$ is a functor in~$C$.  (In super setting
Weil's approach called {\it the language of points} or was
rediscovered in~\cite{L1}
as {\it families}, see also \cite{M,D}.)  This definition might look terribly complicated, but
fortunately one can show that the correspondence
${\mathcal L}\longleftrightarrow L$ is one-to-one and the Lie algebra ${\mathcal L}(C)$,
also denoted $L(C)$, admits a very simple description: $L(C)=(L\otimes
C)_{\bar 0}$.

A {\it Lie superalgebra homomorphism} $\rho: L_1 \longrightarrow L_2$ in these
terms is a functor morphism, i.e., a collection of Lie algebra
homomorphisms $\rho_C: L_1 (C)\longrightarrow L_2(C)$ compatible with morphisms
of supercommutative superalgebras $C\longrightarrow C'$.  In particular, a {\it
representation} of a Lie superalgebra $L$ in a superspace $V$ is a
homomorphism $\rho: L\longrightarrow {\mathfrak{gl}} (V)$, i.e., a~collection of Lie algebra
homomorphisms $\rho_C: L(C) \longrightarrow ( {\mathfrak{gl}} (V )\otimes C)_{\bar 0}$.

\begin{example} Consider a representation
$\rho:{\mathfrak{g}}\longrightarrow{\mathfrak{gl}}(V)$.  The tangent
space of the moduli superspace of deformations of $\rho$ is isomorphic
to $H^1({\mathfrak{g}}; V\otimes V^*)$.  For example, if
${\mathfrak{g}}$ is the $0|n$-dimensional (i.e., purely odd) Lie
superalgebra (with the only bracket possible: identically equal to
zero), its only irreducible representations are the trivial one, {\bf
1}, and $\Pi({\bf 1})$.  Clearly, ${\bf 1}\otimes {\bf 1}^*\simeq
\Pi({\bf 1})\otimes \Pi({\bf 1})^*\simeq {\bf 1}$, and because the
superalgebra is commutative, the differential in the cochain complex
is trivial.  Therefore, $H^1({\mathfrak{g}}; {\bf
1})=\Lambda^1({\mathfrak{g}}^*)\simeq{\mathfrak{g}}^*$, so there are
$\dim {\mathfrak{g}}$ odd parameters of deformations of the trivial
representation.  If we consider~${\mathfrak{g}}$ ``naively'' all of
the odd parameters will be lost.

Which of these infinitesimal deformations can be extended to a global
one is a separate much tougher question, usually solved {\it ad hoc}.

Thus, ${\mathfrak q}tr$ is not a representation of ${\mathfrak{q}}(n)$ according
to the naive definition (``a representation is a Lie superalgebra
homomorphism'', hence, an even map), but is a representation,
moreover, an irreducible one, if we consider odd parameters.
\end{example}

\subsection*{Acknowledgements}
We are thankful to D~Alekseevsky, J~Bernstein, A~Goncharov,
A~Onishchik and I~Shche\-poch\-kina for help; special thanks are
due to P~Grozman who checked computations for the exceptional
superspaces. During preparation of the manuscript we were partly
supported by SFB-170 (1990); DL was also supported by I~Bendixson
grant (1987); NFR, Sweden, and NSF grant DMS-8610730; EP was
supported by the Anna-Greta and Holder Crafoords fond, The Royal
Swedish Academy of Sciences. The final version was moulded during
DL's stay at MPIM-Bonn, to which all of us, and Grozman, are
thankful.

\label{leites-lastpage}
\end{document}